\documentclass[12pt]{amsart}
\usepackage{amssymb}

  \newcommand{\n}{{\bf n}}
  \newcommand{\B}{\mathcal{B}}
 \newcommand{\cB}{\mathcal{X}}
\newcommand{\cX}{\mathcal{X}}
\newcommand{\cY}{\mathcal{Y}}

\newcommand{\cK}{\mathcal{K}}
\newcommand{\G}{\mathcal{G}}
\newcommand{\cZ}{\mathcal{Z}}
\newcommand{\X}{\mathcal{X}}
\newcommand{\gG}{\Gamma}
\newcommand{\C}{\mathbb{C}}
\newcommand{\T}{\mathbb{T}}
\newcommand{\R}{\mathbb{R}}
\newcommand{\E}{\mathbb{E}}
\newcommand{\N}{\mathbb{N}}
\newcommand{\bbZ}{\mathbb{Z}}
\newcommand{\Z}{\mathbb{Z}}
\newcommand{\norm}[1]{\left\Vert #1\right\Vert}
\newcommand{\nnorm}[1]{\lvert\!|\!| #1|\!|\!\rvert}

\theoremstyle{plain}
\newtheorem{theorem}{Theorem}[section]
\newtheorem{lemma}[theorem]{Lemma}
\newtheorem{proposition}[theorem]{Proposition}

\newtheorem*{theorem*}{Theorem}
\newtheorem*{Correspondence1}{Furstenberg's Correspondence Principle}

\newtheorem*{corollary*}{Corollary}
\newtheorem{corollary}[theorem]{Corollary}
\newtheorem*{theoremA}{Theorem A}
\newtheorem*{theoremB}{Theorem B}
\newtheorem*{theoremC}{Theorem C}
\newtheorem*{theoremC'}{Theorem C'}
\newtheorem*{theoremD}{Theorem D}
\newtheorem*{theoremE}{Theorem E}
\newtheorem*{theoremF}{Theorem F}

\theoremstyle{definition}
\newtheorem{definition}[theorem]{Definition}
\newtheorem{example}{Example}

\theoremstyle{remark}

\begin{document}
\title{Multiple ergodic averages for three polynomials and applications}
\author{Nikos Frantzikinakis}
\address[]{Department of Mathematics,
University of Memphis, Memphis, TN 38152-3240}
\email{frantzikinakis@gmail.com}

\begin{abstract}
We find the smallest characteristic factor and a limit formula for
the multiple ergodic averages associated to  any family of three
polynomials and  polynomial families of the form
$\{l_1p,l_2p,\ldots,l_kp\}$. We then derive several multiple
recurrence results and  combinatorial implications, including an
answer to a question of Brown, Graham, and Landman, and a
generalization of the Polynomial Szemer\'edi Theorem of Bergelson
and Leibman for families of three polynomials with not necessarily
zero constant term. We also simplify and generalize  a recent
result of Bergelson, Host, and Kra, showing that for all
$\varepsilon>0$ and every subset of the integers $\Lambda$ the set
$$
\big\{n\in\N\colon d^*\big(\Lambda\cap (\Lambda+p_1(n))\cap
(\Lambda+p_2(n))\cap (\Lambda+
p_3(n))\big)>(d^*(\Lambda))^4-\varepsilon\big\}
$$
has bounded gaps for ``most" choices of integer polynomials
$p_1,p_2,p_3$.
\end{abstract}
\thanks{The  author was partially supported
 by NSF grant DMS-0111298.}



\maketitle

\setcounter{tocdepth}{1}

\tableofcontents

\section{Introduction and main results}
\subsection{Background}
A far reaching generalization of  the  theorem of
Szemer\'edi~\cite{Sz} on arithmetic progressions states that every
subset of the integers with positive upper Banach
density\footnote{If $\Lambda\subset \N$ we define the {\em upper
Banach density} of $\Lambda $ to be
$d^*(\Lambda)=\lim_{N\to\infty}\sup_{M\in\mathbb{N}}|\Lambda\cap[M,M+N)|/N$.}
 contains infinitely many configurations of
the form $\{x, x+p_1(n), \ldots, x+p_k(n)\}$, where
$p_1,\ldots,p_k$ is any collection of \emph{integer polynomials}
(meaning they have integer coefficients) with zero constant term.
This was proved  by Bergelson and Leibman~\cite{BL} using a
Correspondence Principle of Furstenberg and the following result
in ergodic theory:
\begin{theorem}[{\bf Bergelson \& Leibman~\cite{BL}}]\label{T:PSzemeredi}
Let $(X,\mathcal{X},\mu,T)$ be  an invertible measure preserving
system  and let $p_1,\dots, p_k$ be integer polynomials with
$p_i(0) = 0$ for $i = 1, \ldots, k$. If $A\in\B$ with $\mu(A)>0$,
then
\begin{equation}
\label{E:polysl}
\liminf_{N-M\to\infty}\frac{1}{N-M}\sum_{n=M}^{N-1}\mu\bigl(A \cap
T^{p_1(n)}A \cap\ldots\cap T^{p_k(n)}A \bigr) > 0.
\end{equation}
\end{theorem}
A key step in establishing multiple recurrence properties like the
one above is to analyze  the limiting behavior of some closely
related multiple ergodic averages. For the previous result the
relevant ones are
\begin{equation*}
 (P) \quad\qquad\qquad\qquad \qquad\qquad \frac{1}{N-M} \sum_{n=M}^{N-1}
 \ T^{p_1(n)}f_1 \cdot \ldots \cdot  T^{p_k(n)}f_k. \qquad \qquad
 \qquad\qquad\qquad\quad
\end{equation*}
Bergelson and Leibman studied these averages in \cite{BL}, in a
depth that was sufficient for proving \eqref{E:polysl}. Obtaining
a better understanding of their  limiting behavior (as
$N-M\to\infty$) in $L^2(\mu)$ has been a driving force of research
in ergodic theory during the last two decades. The basic approach
for studying them goes back to the original paper of Furstenberg
\cite{Fu1}. Using modern terminology, it consists of finding an
appropriate  factor $\mathcal{C}$ of a given system, called
\emph{characteristic factor}, such that the $L^2$-limit of the
averages in question remains unchanged when each function is
replaced by its projection on this factor. Equivalently, this
means that the averages $(P)$ converge to $0$ in $L^2(\mu)$ as
$N-M\to\infty$ whenever $\E(f_i|\mathcal{C})=0$ for some
$i=1,\ldots,k$, where $\E(f|\mathcal{C})$ is the conditional
expectation of $f$ given $\mathcal{C}$. The next step is to obtain
a concrete description for some well chosen characteristic factor
that is going to facilitate our study.
Using methods from \cite{HK}, this was done in \cite{HK2} for weak
convergence, and in \cite{L2} for strong convergence of the
averages $(P)$.
\begin{theorem}[{\bf Host \& Kra~\cite{HK2}-Leibman~\cite{L2}}]\label{T:HKL}
Let  $p_1,\ldots,p_k$ be a family of essentially distinct
$($meaning, $p_i$ and $p_i-p_j\neq\text{const}$ for $i\neq j)$
integer polynomials. Then there exists a
$d=d(p_1,\ldots,p_k)\in\N$ with the following property: For every
invertible ergodic system some characteristic factor for the
averages $(P)$ is  an inverse limit of $d$-step nilsystems
$($defined in Section~\ref{S:background}$)$.
\end{theorem}

This result opens up the road for a better understanding of the
limiting behavior of the averages $(P)$, and in fact combined with
a recent result of Leibman~\cite{L}  immediately implies that they
converge in $L^2(\mu)$. But  we are still left with some
interesting problems  since computing the smallest characteristic
factor and the actual limit in the case of a nilsystem is still a
difficult task. For example, it is not even clear from the results
in \cite{HK2} and \cite{L2} whether the minimal $d(p_1,p_2)$ is
bounded when the polynomials $p_1,p_2$ vary, and what the limit of
the averages $(P)$  is for $k=2$. Formulas for the limit are known
 when all the polynomials are
linear (see \cite{Z}) or linearly independent (see
\cite{FK}). Also,  very recently some other cases where covered in  \cite{L4}.

In this paper we are going  to
find the smallest characteristic factor and limit formulas for the
averages $(P)$ for any family of three polynomials and for
polynomial families of the form $\{l_1p,l_2p,\ldots,l_kp\}$.
We will then use these results to
derive several combinatorial implications.





\subsection{Results in ergodic theory}
Given a measure preserving system and a family of integer
polynomials $P=\{p_1,\ldots,p_k\}$ we say that a factor
$\mathcal{C}$ is the \emph{smallest characteristic factor for
$P$}, if it is a characteristic factor for the averages $(P)$ and
it is a factor of every other such characteristic factor. We will
completely determine the structure of the smallest characteristic
factor for any family of three polynomials and the family
$\{l_1p,l_2p,\ldots,l_kp\}$. The reader who is not familiar with
the notions we use in ergodic theory may wish to consult
Section~\ref{SS:ergodic} first.

We first deal with  the polynomial family
$\{l_1p,l_2p,\ldots,l_kp\}$:
\begin{theoremA}\label{T:p(n)}
Let $(X,\mathcal{X},\mu,T)$ be an invertible  ergodic system, $p$
be a nonconstant integer polynomial, and $l_1,\ldots,l_k$ nonzero
distinct integers. If $k\geq 2$ then the $(k-1)$-step nilfactor
$\mathcal{Z}_{k-1}$ is the smallest characteristic factor for the
multiple ergodic averages
$$
\frac{1}{N-M}\sum_{n=M}^{N-1} T^{l_1p(n)}f_1\cdot
T^{l_2p(n)}f_2\cdot\ldots\cdot T^{l_kp(n)}f_k. \footnote{If $k=1$
it is well known (\cite{Fu2}) that the rational Kronecker factor
$\mathcal{K}_{rat}$ is a characteristic factor.}
$$
 Moreover, if the system is totally ergodic
then the $L^2$-limit as $N-M\to\infty$ does not depend on the
choice of the polynomial $p$ and can be computed explicitly.
\end{theoremA}
We will use this result to answer a question of Brown, Graham and
Landman~\cite{BGL} (see Theorem~D), and to deal with
characteristic factors for families of three polynomials (see
Theorem~B). The proof of Theorem~A is based on
Proposition~\ref{L:p(n)} which enables us to compare the family
$\{l_1p,l_2p,\ldots,l_kp\}$ with the family
$\{l_1n,l_2n,\ldots,l_kn\}$.

Before  we deal with  families of three polynomials we take a
moment to define three classes of polynomial families that will
help us expedite the  discussion:
\begin{definition}
We say that the family $\{p_1,p_2,p_3\}$  of essentially distinct
integer polynomials is of type $(e_1)$, $(e_2)$, $(e_3)$, if some
permutation of the polynomials
$\{\tilde{p}_1,\tilde{p}_2,\tilde{p}_3\}$, where
$\tilde{p}_i=p_i-p_i(0)$, $i=1,2,3$, has the form  $\{lp,mp,rp\}$,
$\{lp,mp,kp^2+rp\}$, $\{kp^2+lp,kp^2+mp,kp^2+rp\}$
correspondingly, for some integer polynomial $p$ and constants
$k,l,m,r\in\Z$ with $k\neq 0$.
\end{definition}
\begin{theoremB}\label{T:3polys}
Let $(X,\mathcal{B},\mu,T)$ be an invertible ergodic  system and
$\{p_1,p_2,p_3\}$ be a family of essentially distinct integer
polynomials. Consider the multiple ergodic averages
\begin{equation}\label{E:3polys}
\frac{1}{N-M}\sum_{n=M}^{N-1} T^{p_1(n)}f_1\cdot
T^{p_2(n)}f_2\cdot T^{p_3(n)}f_3.
\end{equation}
Then  the following three mutually exclusive cases describe the
smallest\footnote{In case $(i)$ it is the smallest under the extra
assumption that the system is totally ergodic.} characteristic
factor for the averages $\eqref{E:3polys}:$

\noindent $(i)$ It is the rational Kronecker factor
$\mathcal{K}_{rat}$ if the polynomials
$\tilde{p}_1,\tilde{p}_2,\tilde{p}_3$ are linearly independent.

\noindent $(ii)$ It is the $2$-step nilfactor $\cZ_2$ if the
polynomials are of type $(e_1)$, and the $2$-step affine factor
$\mathcal{A}_2$ if  the polynomials are of type $(e_2),(e_3)$.

\noindent $(iii)$ It is the Kronecker factor $\mathcal{K}$ in all
other cases.

\noindent Furthermore, if the system  is totally ergodic we can
give explicit formulas for the $L^2$-limit of the averages
\eqref{E:3polys} as $N-M\to\infty$.
\end{theoremB}

The following examples illustrate the different  limiting
behaviors the averages \eqref{E:3polys} may exhibit:

$(a)$ If  $P=\{n,n^2, n^3\}$ then  the rational Kronecker factor
$\mathcal{K}_{rat}$ is characteristic. In the totally ergodic case
the limit is the product of the integrals of the three functions.

$(b)$ If  $P=\{n,n^2,n^2+n\}$ then  the Kronecker factor
$\mathcal{K}$ is characteristic. In the totally ergodic case the
limit is the same as in the case of the double averages (averaging
over  $m,n$) associated to the family $\{m,n,m+n\}$.

$(c)$ If $P=\{n,2n,n^3\}$ then  the Kronecker factor $\mathcal{K}$
is characteristic.  In the totally ergodic case the limit is the
product of  the limit of the  ergodic averages corresponding to
the family  $\{n,2n\}$ and the integral of the third function.

$(d)$ If  $P=\{n, 2n, n^2\}$ then  the two step affine factor
$\mathcal{A}_2$ is characteristic. This is the first example that
we know of a polynomial family  with smallest characteristic
factor (for totally ergodic systems)  not of the form
$\mathcal{Z}_m$ for some nonnegative integer $m$. In the totally
ergodic case the limit can be computed explicitly and unlike the
case $\{n, 2n, n^3\}$ it depends nontrivially on the third
function.

$(e)$ If $P=\{n,2n, 3n\}$ or $P=\{n^2,2n^2, 3n^2\}$ then the
$2$-step nilfactor $\mathcal{Z}_2$ is characteristic. In the
totally ergodic case the limit is the same in both cases and can
be computed explicitly.

The proof of Theorem~B is rather complicated so let us briefly
explain the main ideas. Our first step is  to use
Theorem~\ref{T:HKL}  in order to show that it suffices to restrict
our study to totally ergodic nilsystems (Proposition~\ref{P:te}).
At this point we are left with establishing various  uniform
distribution properties on nilmanifolds. Our main
 tool  is a ``\emph{reduction to affine argument}" which consists of the
 following two
steps: $(i)$ Reduce the uniform distribution problem to a simpler
one that involves only nilpotent affine transformations on finite
dimensional tori. This reduction is done using Theorems~\ref{T:L}
and \ref{T:FK} but varies in difficulty depending on the problem.
$(ii)$ Verify the simplified (but often challenging) problem ``by
hand" for affine transformations. Here is a more detailed sketch
of how this plan is executed to deal with the various parts of
Theorem~B:

  Part $(i)$ deals with
linearly independent polynomial families, a case that has been
already worked out  in~\cite{FK} and \cite{FK2} using a
``reduction to affine argument''. Part $(ii)$ deals with families
of type $(e_1),(e_2),$ and $(e_3)$. A typical family of type
$(e_1)$ is $P=\{n^2,2n^2,3n^2\}$. It follows from Theorem~A (which
is again proved using a ``reduction to affine argument'') that
$P\sim\{n,2n,3n\}$.\footnote{We say that two polynomial families
are equivalent (we write $P\sim Q$) if for totally ergodic systems
the corresponding averages \eqref{E:3polys} have the same
$L^2$-limit.} This completes our task since it is known
(\cite{CL1}, \cite{CL2}) that for this family the factor
$\mathcal{Z}_2$ is characteristic. A typical family of type
$(e_2)$ is  $P=\{n,2n,n^2\}$. To deal with this case, we first use
Van der Corput Lemma and the fact that for families of the form
$\{n,n^2,n^2+kn\}$ ($k\neq 0$) the Kronecker factor $\mathcal{K}$
is characteristic (Lemma~\ref{L:w2'} and \ref{L:w2}), to show that
if $\E(f_3|\mathcal{K})=0$ then the averages \eqref{E:3polys}
converge to zero in $L^2$. This fact greatly simplifies the
analysis,  and we are led  to consider averages corresponding to
the family $\{n,2n\}$ for a transformation $S=T\times R$ where $R$
is a $2$-step  affine transformation on $\T^2$. From this we
deduce using a result  from \cite{Fr} that the factor
$\mathcal{A}_2$ is characteristic. To complete the study of
families of type $(e_2)$ we need also to show that
$\{lp,mp,kp^2+rp\}\sim \{ln,mn,kn^2+rn\}$ when the polynomial  $p$
is nonconstant. To do this we use Proposition~\ref{L:p(n)} (again
proved using a``reduction to affine argument'') which roughly
speaking tells us that if $p(n)$ is a nonconstant integer
polynomial, then  the substitution $n\to p(n)$ does not change the
distribution of any polynomial sequence that has connected orbit
closure.
Families of type $(e_3)$ are easily reduced to families of type
$(e_2)$, thus completing the study of part $(ii)$. Finally, to
deal with part $(iii)$, the crucial step is
Proposition~\ref{T:weyl}. We show there  that the polynomial
families that were not covered by part $(i)$ and $(ii)$ have Weyl
complexity $2$ (defined in Section~\ref{S:weyl}). This fact,
combined with Lemmas~\ref{L:w2'} and \ref{L:w2}, allows us to
conclude that in this case the Kronecker factor $\mathcal{K}$ is
characteristic.



The following is an immediate corollary of Theorem~B:
\begin{corollary*}
For any two essentially distinct  polynomials and every invertible
ergodic system, the Kronecker factor $\mathcal{K}$ is
characteristic for the corresponding  averages $(P)$, and for any
three essentially distinct polynomials the $2$-step nilfactor
$\mathcal{Z}_2$ is characteristic.
\end{corollary*}
It seems plausible that for  $k\geq 2$ the $(k-1)$-step nilfactor
$\mathcal{Z}_{k-1}$ is characteristic for any family of $k$
essentially distinct  polynomials. Moreover, one would expect that
for $k\geq 2$ the smallest $m$ for which the factor $\cZ_{m-1}$ is
characteristic for a family $P$ of essentially distinct integer
polynomials is $W(P)$ (defined in Section~\ref{S:weyl}). It is an
immediate consequence of Theorem~B and Proposition~\ref{T:weyl}
that both statements hold for $k=2,3$.

Next we establish a multiple recurrence result that generalizes a
result of Bergelson, Host, and Kra~\cite{BHK}:

\begin{theoremC}\label{T:lower1} Let $(X,\mathcal{X},\mu,T)$ be
an invertible ergodic  system, $A\in\mathcal{X}$ with $\mu(A)>0$,
and $\{p_1,p_2,p_3\}$ be  integer polynomials with $p_i(0)=0$, for
$i=1,2,3$. Then for every $\varepsilon>0$ the set
$$
\left\{n\in\N\colon \quad \mu\left( A\cap T^{p_1(n)}A\cap
T^{p_2(n)}A\right)\geq \mu(A)^{3}-\varepsilon\right\}
$$
has bounded gaps. Moreover, the set
$$
\left\{n\in\N\colon \quad \mu\left( A\cap T^{p_1(n)}A\cap
T^{p_2(n)}A\cap T^{p_3(n)}A\right)\geq
\mu(A)^{4}-\varepsilon\right\}
$$
has bounded gaps,
 unless   the polynomials are essentially distinct and of type $(e_1)$
 with $l<m<r$ and $r\neq l+m$,  or of type $(e_2)$, $(e_3)$.
\end{theoremC}
  This result was established in \cite{BHK} for the polynomial
 families $\{n,2n\}$ and $\{n,2n,3n\}$. Moreover, it was shown
 that an analogous result fails for the
family $\{n,2n,3n,4n\}$, in fact no fixed power of $\mu(A)$ works
as a lower bound. To prove Theorem~C
we use Theorem~A and parts $(i), (iii)$ of Theorem~B. We remark
that even for the two cases covered
 in \cite{BHK} our argument is different and much simpler ($1$ and
 $2$ pages long correspondingly). The crucial  observation is that
although we cannot get good lower bounds for the averages
corresponding to the families $\{n,2n\}$ and $\{n,2n,3n\}$ if we
average over the full set of positive integers, we can get optimal
lower bounds as long as the average is taken  over an
appropriately chosen subset of the integers (that depends on the
system given).\footnote{This is best exemplified by considering an
irrational rotation $\alpha\in\T$. Although it is not possible to
get good lower bounds for the average of the sequence $\mu(A\cap
(A+n\alpha)\cap (A+2n\alpha))$ when $n$ ranges over $\N$, if we
restrict the range of $n$ to the set
$S_\delta=\{n\colon\{n\alpha\}\leq \delta\}$, it is easy to show
that  for every $\varepsilon>0$ if $\delta$ is small enough then
the average of the corresponding subsequence   is greater than
$\mu(A)^3-\varepsilon$.} This observation greatly simplifies the
whole analysis, as we do not have to rely on the rather
complicated nilsequence decompositions used in \cite{BHK}.

For the exceptional polynomial families of Theorem~C  we believe
that the analogous result fails
 and we provide conditional
counterexamples in Section~\ref{counterexamples}.
\subsection{Results in combinatorics}\label{S:comb_results}
We are going to utilize the previous results in ergodic theory to
derive several  implications in combinatorics. We mention them in
increasing degree of difficulty.

We start with an answer to a question of Brown, Graham, and
Landman. In \cite{BGL} the authors define a set $S\subset \Z$ to
be large if every finite coloring of the positive integers
contains arbitrarily long monochromatic arithmetic progressions
with common difference a nonzero integer in $S$. It follows from
Theorem~\ref{T:PSzemeredi} that if $p$ is an integer polynomial
with $p(0)=0$ then the set $S_p=\{p(n)\colon n\in \N\}$ is large.
If we do not assume that $p(0)= 0$ an obvious necessary condition
for the set $S_p$ to be large is that it contains multiples of
every positive integer. The authors of \cite{BGL} asked whether
this condition is also sufficient and  in particular whether the
range of the polynomial $p(n)=(n^2-13)(n^2-17)(n^2-221)$ is large;
this is an example of a polynomial with no linear integer factors
whose range does contain multiples of every positive
integer\footnote{As shown in \cite{Be}, the smallest possible
degree of a polynomial having this property is $5$, an example is
$p(n)=(n^3-19)(n^2+n+1)$.} (this can be easily verified using
properties of the Legendre symbol). We will give a positive answer
to these questions,  in fact we will verify a stronger ``density"
statement. We say that $S\subset \Z$ is a \emph{set of multiple
recurrence} if every subset of the integers with positive density
contains arbitrarily long arithmetic progressions with common
difference a nonzero integer in $S$. We show:
\begin{theoremD}\label{T:polys}
Let $p$ be an integer polynomial. Then $S_p=\{p(n)\colon n\in\N\}$
is a set  of multiple recurrence if and only if it contains
multiples of every positive integer.
\end{theoremD}
To prove this result we use Theorem~A and  Furstenberg's Multiple
Recurrence Theorem~\cite{Fu1}. Polynomials that satisfy the
conditions of Theorem~D have been studied in \cite{Be}. It is
shown there that  $p(n)\equiv 0\pmod m$ is solvable for every
$m\in\N$ if and only if it is solvable for a finite set of
$m\in\N$ explicitly depending on $p$.

Our next application is to construct a set $S$ that has bad
recurrence properties but its set of squares $S^2$  is a set of
multiple recurrence.  Note that if $S$ is a set multiple recurrence
it is not known whether its set of squares $S^2$ is always a set of
multiple recurrence
 (the chromatic
version of this question was conjectured to be true in
\cite{BGL}).
\begin{theoremE}\label{T:goodbad}
There exists a set $S\subset \N$ that is not a set of multiple (in
fact not even  single) recurrence but $p(S)=\{p(s),s\in S\}$ is a
set of multiple recurrence for every integer polynomial $p$ with
degree greater than $1$.
\end{theoremE}
Our example is explicit, in fact we show that the set
$S=\big\{n\in\N\colon \{n \sqrt{2}\}\in [1/4,3/4]\big\}$ works. To
prove  this we rely on Lemma~\ref{L:squares}.

Our next application deals with an extension of
Theorem~\ref{T:PSzemeredi} to families of polynomials with not
necessarily zero constant term. We say that the family of integer
polynomials $\{p_1,\ldots,p_k\}$ is \emph{universal} if every
subset of the integers with positive density contains infinitely
many configurations of the form $\{x,x+p_1(n),\ldots,x+p_k(n)\}$,
where $x,n\in \N$. From Theorem~\ref{T:PSzemeredi} we know that
every family of integer polynomials with zero constant term is
universal. We show:

\begin{theoremF}\label{T:universal}
The family of integer polynomials $\{p_1,p_2,p_3\}$  is universal
if and only if the congruence $p_1(n)\equiv p_2(n)\equiv
p_3(n)\equiv 0\pmod{m} $ has a solution for every $m\in\N$.
\end{theoremF}
To prove this result we make essential use of Theorem~B, so we are
currently unable to extend it  to deal with families of $k$ polynomials for
$k\geq 4$.

Finally, using a modification of the Correspondence Principle of
Furstenberg, due to Lesigne (see Section~\ref{S:applications}), we
give the following combinatorial implication of  Theorem~C:
\begin{theoremC'}\label{T:lower2} Let
$\Lambda\subset \N$ and $p_1,p_2,p_3$ be  integer polynomials with
$p_i(0)=0$ for $i=1,2,3$. Then for every $\varepsilon>0$ the set
$$
\{n\in\N\colon \quad
d^*\bigl(\Lambda\cap(\Lambda+p_1(n))\cap(\Lambda
+p_2(n))\bigr)\geq (d^*(\Lambda))^{3}-\varepsilon\}
$$
has bounded gaps, and the set
$$
\{n\in\N\colon \quad
d^*\bigl(\Lambda\cap(\Lambda+p_1(n))\cap(\Lambda
+p_2(n))\cap(\Lambda +p_3(n))\bigr)\geq
(d^*(\Lambda))^{4}-\varepsilon\}
$$
has bounded gaps,  unless   the polynomials are essentially
distinct and of type $(e_1)$
 with $l<m<r$ and $r\neq l+m$,  or of type $(e_2)$, $(e_3)$.
\end{theoremC'}
  Examples of random sets show that the lower bounds given are tight.
  The same result was established in \cite{BHK} in the
special case of the polynomial families $\{n,2n\}$ and
$\{n,2n,3n\}$. In the case of the family $\{n,2n\}$ a related
finite version of this result was established by Green~\cite{G}.
 Some other examples of eligible $3$-term polynomial
families are the following: $\{n,3n,4n\}$, $\{n^k,2n^k,3n^k\}$ for
all $k\in \N$, $\{n,n^2,an^2+bn\}$ with $a\neq 0$, and
$\{n,2n,n^k\}$ for all $k\geq 3$. It was shown in \cite{BHK} that
similar lower bounds fail for the polynomial family
$\{n,2n,3n,4n\}$. In  contrast to this, similar lower bounds hold
for any family of $k$ linearly independent polynomials with zero
constant term (see \cite{FK2}).

As was the case with the corresponding result in ergodic theory,
for the exceptional polynomial families of Theorem~C' we believe
that the analogous result fails
and we provide conditional
counterexamples in Section~\ref{counterexamples}.

{\bf Notation:} The following notation will be used throughout the
article: $Tf=f\circ T$, $e(x)=e^{2\pi i x}$, $\{x\}=x-[x]$,
$\text{UD-}\!\lim(a_n)=0$ if for every $\varepsilon>0$ we have
$d^*(\{n\colon |a_n|>\varepsilon\})=0$.

{\bf Acknowledgements.} The author would like to thank B. Kra for
 helpful discussions during the preparation of this article, M.
 Johnson for helpful remarks,
and S. Leibman for providing the simple proof of
Proposition~\ref{L:p(n)}.

\section{Background in ergodic theory and nilsystems}\label{S:background}

\subsection{Ergodic theory background and notation}\label{SS:ergodic}
Background information we assume in this article can be found in
the books \cite{Fu2}, \cite{Pe}, \cite{Wa}. By a \emph{measure
preserving system} (or just \emph{system}) we mean a quadruple
$(X,\mathcal{X},\mu, T)$, where $(X,\mathcal{X},\mu)$ is a
probability space and $T\colon X\to X$ is a measurable map such
that $\mu(T^{-1}A)=\mu(A)$ for all $A\in\mathcal{X}$. Without loss
of generality we can assume that the probability space is
Lebesgue.
 A \emph{factor}
of a  system  can be defined in any of the following three  ways:
it is a $T$-invariant sub-$\sigma$-algebra $\mathcal{D}$ of $\cX$,
it is a $T$-invariant sub-algebra $\mathcal{F}$ of $L^\infty(X)$,
or it is a system $(Y, \cY, \nu, S)$ and a measurable map
$\pi\colon X'\to Y'$, where $X'$ is a $T$-invariant set and $Y'$
is an $S$-invariant set of full measure, such that
$\mu\circ\pi^{-1} = \nu$ and $S\circ\pi(x) = \pi\circ T(x)$ for
$x\in X'$. . In a slight abuse of terminology, when any of these
conditions holds, we say that $Y$ (or the appropriate
$\sigma$-algebra of $\cX$) is a factor of $X$ and call $\pi\colon
X'\to Y'$ the factor map. If the  factor map $\pi\colon X'\to Y'$
can be chosen to be injective, then we say that the systems
$(X,\cB, \mu, T)$ and $(Y, \cY, \nu, S)$ are \emph{isomorphic}
(bijective maps on Lebesgue spaces have measurable inverses).

If $\cY$ is a $T$-invariant sub-$\sigma$-algebra of $\cX$ and
$f\in L^2(\mu)$, we define the \emph{conditional expectation
$\mathbb{E}(f|\cY)$ of $f$ with respect to $\cY$} to be the
orthogonal projection of $f$ onto $L^2(\cY)$. We frequently use
the identities
$$
\int \mathbb{E}(f|\cY) \ d\mu= \int f\ d\mu, \quad
T\,\mathbb{E}(f|\cY)=\mathbb{E}(Tf|\cY).
$$
For each $r\in\mathbb{N}$, we define $\mathcal{K}_{r}$ to be the
factor induced by the algebra
$$\{f\in L^\infty(\mu):T^rf=f\}\ .$$
We define $\mathcal{K}_{rat}$ to be the factor induced by the
algebra generated by the functions
$$
\{f\in L^\infty(\mu):T^rf=f \text{ for some } r\in\mathbb{N}\}\ .
$$
The Kronecker factor $\mathcal{K}$  is induced by the algebra
spanned by the bounded eigenfunctions of $T$, i.e. functions that
satisfy $Tf=e(a)\cdot f$ for some $a\in \R$.  We also define
higher order eigenfunctions and their corresponding factors. Let
$\mathcal{E}_0$ denote the set of eigenvalues of $T$ and for
$k\in\N$ we define inductively
$$
\mathcal{E}_k=\{f\in L^\infty(\mu)\colon |f|=1 \text{ and
 } Tf \cdot \bar{f}\in\mathcal{E}_{k-1}(T)\}.$$
We call the factor spanned by $\mathcal{E}_k$ the $k$-\emph{step
affine} factor of the system, and denote it by $\mathcal{A}_k$.
The reason for this notation is that for totally ergodic systems
the factor system induced by  $\mathcal{A}_k$ is isomorphic to a
nilpotent $k$-step affine transformation on some connected compact
abelian group (this is a result of Abramov~\cite{A}), and
$\mathcal{A}_k$ is the largest factor with this property.

The transformation $T$ is \emph{ergodic} if $\cK_1$ consists only
of constant functions, and $T$ is \emph{totally ergodic} if
$\cK_{rat}$ consists only of constant functions.
Every system $(X,\cX,\mu,T)$ has an \emph{ergodic decomposition},
meaning that we can write $\mu=\int \mu_t\ d\lambda(t)$, where
$\lambda$ is a probability measure on $[0,1]$ and $\mu_t$ are
$T$-invariant probability measures on $(X,\cX)$ such that the
systems $(X,\cX,\mu_t,T)$ are ergodic for $t\in [0,1]$. We
sometimes denote the ergodic components by $(T_t)_{t\in [0,1]}$.

We say that the system $(X,\cX,\mu,T)$ is an \emph{inverse limit
of a sequence of factors} $(X,\cX_j,\mu,T)$ if
$\{\cX_j\}_{i\in\mathbb{N}}$ is an increasing sequence of
$T$-invariant sub-$\sigma$-algebras such that
$\bigvee_{j\in\N}\mathcal{X}_j=\mathcal{X}$ up to sets of measure
zero.

Following \cite{HK}, for every  system $(X,\mathcal{X},\mu,T)$ and
function $f\in L^\infty(\mu)$, we define inductively the seminorms
$\nnorm{f}_k$ as follows: For $k=1$ we set
$\nnorm{f}_1=|\E(f|\mathcal{I})|$,\footnote{In \cite{HK} the
authors work with ergodic systems, in which case $\nnorm{f}_1=\int
f \ d\mu$, but the whole discussion can be carried out for
nonergodic systems as well without extra difficulties.} where
$\mathcal{I}$ is the $\sigma$-algebra of $T$-invariant sets. For
$k\geq 2$ we set
\begin{equation}
\label{eq:recur} \nnorm f_{k+1}^{2^{k+1}} =\lim_{N\to+\infty}\frac
1N\sum_{n=0}^{N-1} \nnorm{\overline{f}\cdot T^nf}_{k}^{2^{k}}.
\end{equation}
It was shown in~\cite{HK} that for every integer $k\geq 1$,
$\nnorm\cdot_k$ is a seminorm on $L^\infty(\mu)$ and  it defines
factors $\cZ_{k-1}$ in the following manner: the $T$-invariant
sub-$\sigma$-algebra $\cZ_{k-1}$ is characterized by
$$
\text{ for } f\in L^\infty(\mu),\  \E(f|\cZ_{k-1})=0\text{ if and
only if } \nnorm f_{k} = 0 .
$$
We remark that if $(T_t)_{t\in [0,1]}$ are the ergodic components
of the system then $\E(f|\mathcal{Z}_k(T))=0$ if and only if
$\E(f|\mathcal{Z}_k(T_t))=0$ for a.e. $t\in[0,1]$. For ergodic
systems the factor  $\cZ_0$ is trivial,
$\cZ_1=\mathcal{A}_1=\mathcal{K}$, and
$\mathcal{A}_k\subset\mathcal{Z}_k$ (the inclusion is in general
proper for  $k\geq 2$). The factors $\cZ_k$   are of particular
interest since they are characteristic for $L^2$-convergence of
ergodic averages $(P)$. Moreover, in \cite{HK} it was shown that
the factor $\cZ_k$ is an inverse limit of $k$-step nilsystems
which brings us to our next topic of discussion.

\subsection{Nilsystems, definition and examples}\label{SS:nilmanifolds}
Fundamental properties of nilsystems, related to our discussion,
were studied in \cite{AGH}, \cite{Pa}, \cite{Les}, \cite{L0}, and
\cite{L1}. Below we summarize some facts that we shall  use, all
the proofs can be found in \cite{L}.

 Given a topological group $G$, we denote the identity
element by $e$ and
  we let $G_0$ denote the connected component of $e$.
  If $A, B\subset G$, then $[A,B]$ is
defined to be the subgroup $\{[a,b]:a\in A, b\in B\}$ where
$[a,b]=ab a^{-1}b^{-1}$. We define the commutator subgroups
recursively by $G_1=G$ and
 $G_{k+1}=[G, G_{k}]$. A group
$G$ is said to be {\em $k$-step nilpotent} if its $(k+1)$
commutator $G_{k+1}$ is trivial. If $G$ is a $k$-step nilpotent
Lie group and $\gG$ is a discrete cocompact subgroup, then the
compact space $X = G/\gG$ is said to be a {\em $k$-step
nilmanifold}.  The group $G$ acts on $G/\gG$ by left translation
and the translation by a fixed element $a\in G$ is given by
$T_{a}(g\gG) = (ag) \gG$. Let $m$ denote the unique probability
measure on $X$ that is invariant under the action of $G$ by left
translations (called the {\em Haar measure}) and let $\G/\gG$
denote the Borel $\sigma$-algebra of $G/\gG$. Fixing an element
$a\in G$, we call the system $(G/\gG, \G/\gG, m, T_{a})$ a {\em
$k$-step nilsystem} and call the map $T_a$ a {\em nilrotation}.

If $H$ is a closed subgroup of $G$ then $Y=(H\Gamma)/\Gamma\simeq
H/(H\cap \Gamma)$ may not be compact in general (take $X=\R/\Z$
and $H=\{t\sqrt{2}\colon t\in \R\}$), but if $Hx$ is closed in $X$
for some $x\in X$, then it can be shown that $Y$ is compact and
the set   $Hx$ can be given the structure of a nilmanifold. In
particular if $x=g\Gamma$ for some $g\in G$ we have   $Hx\simeq
H/\Delta$ where $\Delta=H\cap g\Gamma g^{-1}$. We call any such
set a {\em sub-nilmanifold} of $X$.

 Examples of nilsystems are rotations on  compact abelian Lie groups,
 and more generally, every  nilpotent affine transformation on a compact
 abelian Lie group is isomorphic to a nilsystem (see Example 1).
But these examples do not cover all the possible nilsystems (see
Example 2).

\begin{example}
\label{ex:one} On the space $G=\bbZ\times\R^2$,
define multiplication as follows: \\
if $g_1=(m_1,x_1,x_2)$ and $g_2=(n_1,y_1,y_2)$, let
$$
g_1\cdot g_2=(m_1+n_1,x_1+y_1, x_2+y_2+m_1y_1).
$$ Then $G$ is a $2$-step
nilpotent group and  the discrete subgroup $\gG=\bbZ^3$ is
cocompact. If  $a=(m_1,a_1,a_2)$, it turns out that $T_a$ is
isomorphic to the a nilpotent affine transformation $S \colon
\mathbb{T}^2\to\mathbb{T}^2$ given by
$$
S(x_1,x_2)=(x_1+a_1,x_2+m_1x_1+a_2).
$$
\end{example}

\begin{example}
\label{ex:two} On the space $G=\R^3$,
define multiplication as follows: \\
if $g_1=(x_1,x_2,x_3)$ and $g_2=(y_1,y_2,y_3)$, let
$$
g_1\cdot g_2=(x_1+y_1,x_2+y_2, x_3+y_3+x_1y_2).
$$ Then $G$ is a $2$-step
nilpotent group and  the discrete subgroup $\gG=\bbZ^3$ is
cocompact. Let $a=(a_1,a_2,0)$, where $a_1,a_2\in [0,1)$ are
linearly independent.  It turns out that $T_a$ is  isomorphic to a
skew product transformation $S \colon \mathbb{T}^3\to\mathbb{T}^3$
that has  the form
$$
S(x_1,x_2)=(x_1+a_1,x_2+a_2, x_3+f(x_1,x_2)),
$$
where  $f\colon \mathbb{T}^2\to\mathbb{T}$ is defined by
$$
f(x_1,x_2)= (x_1+a_1) [x_2 +a_2] - x_1 [x_2] -a_1x_2.
$$
It can be shown that  $S$ (or $T_a$) is \emph{not} isomorphic to a
nilpotent affine transformation on some finite dimensional torus.
\end{example}

Let $(X=G/\Gamma,\mathcal{G}/\Gamma,m,T_a)$ be an ergodic
nilsystem. The subgroup  $<G_0,a>$ projects to an open subgroup of
$X$ that is invariant under $a$. By ergodicity this projection
equals $X$. Hence, $X=<G_0,a>/\Gamma'$ where $\Gamma'=\Gamma\cap
<G_0,a>$. Using this representation of $X$  for  ergodic
nilsystems we have that
\begin{equation}\label{H}
 \quad  G \emph{ is generated by the connected component of the
identity element  and } a.
\end{equation}
>From now on when we work with an ergodic nilsystem we will freely
assume that hypothesis \eqref{H} is satisfied. We remark that
under this hypothesis it was shown in \cite{L0} that for every
integer $k\geq 2$ the group $G_k$ is connected.

We will make frequent use of the following simple facts:
\begin{proposition}\label{basic}
Let $(X=G/\Gamma,\mathcal{G}/\Gamma,m, T_a)$ be an ergodic
nilsystem. Then

(i) The system is totally ergodic if and only if  $X$ is
connected.

(ii) There exists an $r\in \N$ such that the $($finitely many$)$
ergodic components of $T_a^r$ are totally ergodic.
\end{proposition}
\begin{proof} We first prove statement $(i)$.
Suppose that the system is totally ergodic. Let $X_0$ be the
identity component of $X$. Since $X$ is compact, it is a disjoint
union of a finite number of translations of $X_0$. Since $a$
permutes these copies, there exists an $r\in\N$ such that $a^r$
preserves $X_0$. By assumption the translation by $T_{a^r}=T_a^r$
is ergodic and so $X_0=X$.

Conversely, suppose that $X$ is connected and let $r\in\N$.
Because $T_a$ is ergodic, there exists  $x_0\in X$ such that  the
sequence $\{a^{n}\pi(x_0)\}_{n\in\N}$ is uniformly distributed in
$Z=G/([G,G]\Gamma)$, where $\pi\colon X\to Z$ is the natural
projection. Since $Z$ is a connected compact abelian group, it is
well known that  $\{a^{rn}\pi(x_0)\}_{n\in\N}$ is also uniformly
distributed in $Z$. By  Theorem~\ref{T:L} below
 we have that $T_a^r=T_{a^r}$ is
ergodic. Since $r\in\N$ was arbitrary, $T_a$ is totally ergodic.

We now prove statement $(ii)$. The Kronecker factor of an ergodic
nilsystem is isomorphic to a rotation on a monothetic compact
abelian Lie group $G$. Every such group has the form $\Z_{d_1}
\times \T^{d_2}$ for some positive integer $d_1$ and nonnegative
integer $d_2$, where $\bbZ_d$ denotes the cyclic group with $d$
elements. It follows that $\mathcal{K}_{rat}=\mathcal{K}_{d_1}$,
and $T^{d_1}$ has finitely many ergodic components and they are
all totally ergodic.
\end{proof}

\subsection{Factors of nilsystems}
Given an ergodic  nilsystem the following result allows us to
identify its factors $\mathcal{Z}_k(T_a)$:
\begin{theorem}[{\bf Ziegler \cite{Z1}}] \label{T:Kronecker}
Let  $(X=G/\Gamma, \mathcal{G}/\Gamma,m,T_a)$ be an ergodic
nilsystem.
Then $f\in \mathcal{Z}_k(T_a)$ if and only if $f\in L^{\infty}(m)$
and $f$ factors through $G/(G_{k+1}\Gamma)$.
\end{theorem}
It will also  be convenient for us to identify the $2$-step affine
factor $\mathcal{A}_2$ of an ergodic nilsystem. We adapt a
technique from \cite{L0} to do this. We first need a lemma:

\begin{lemma} \label{L:G_k}
Let  $(X=G/\Gamma, \mathcal{G}/\Gamma,m,T_a)$ be an ergodic
nilsystem. If $f\in\mathcal{E}_2(T_a)$ then for every $b\in G$ we
have $f_b\in \mathcal{E}_2(T_a)$, where $f_b(x)=f(bx)$.
\end{lemma}
\begin{proof}
We know that  $\mathcal{A}_2(T_a)$ is a factor of
$\mathcal{Z}_2(T_a)$, so by Theorem~\ref{T:Kronecker} the function
$f$ factors through $G/(G_3\Gamma)$. Hence,  after replacing $G$
by $G/G_3$ we can assume that $G$ is $2$-step nilpotent. We know
from \cite{A} that $|f|=\text{const}$, so we can assume that
$|f|=1$ in which case we have that $\bar{f}=f^{-1}$. Since
$f\in\mathcal{E}_2(T_a)$ there exists $h\in\mathcal{E}_1(T_a)$
such that
$$
f(ax)=h(x)\cdot f(x).
$$
By Theorem~\ref{T:Kronecker} the function $h$ factors through the
compact abelian group $G/([G,G]\Gamma)$. Moreover, since $h$ is an
eigenfunction of $T_a$ it is a character of $G$.

We first claim that
\begin{equation}\label{E:f_c}
\text{ if }  c\in[G,G] \text{ then } f_c=\lambda_c \cdot f
\end{equation}
for some constant $\lambda_c\in\C$.
 Since $h(cx)=h(x)$ and  $c$
belongs to the center of $G$ we find that
$$
T_af_c=f(cax)=f(acx)=h(cx)\cdot f(cx)=h(x)\cdot f(cx)=h\cdot f_c.
$$
 Hence, $f_c\cdot\bar{f}\in \mathcal{E}_1(T_a)$. We define a map
$\phi\colon [G,G]\to \mathcal{E}_1(T_a)$ by
$\phi(c)=f_c\cdot\bar{f}$. It suffices to show  that
$\phi([G,G])\subset C$, where $C$ is the set of constant
functions. We will use a connectedness argument to do this. If we
equip $\mathcal{E}_1(T_a)$ with the $L^2(m)$ topology then the map
$\phi$ is continuous. Since $T_a$ is ergodic the connected
component of the function $1$ in $\mathcal{E}_1(T_a)$ is the set
$C$. Since $\phi$ is continuous, $\phi(e)=1$, and $[G,G]$ is
connected,  we have that $\phi([G,G])\subset C$. This proves the
claim.

We now show that for every $b\in G$ we have $f_b\in
\mathcal{E}_2(T_a)$. We compute
\begin{equation}\label{E:f1}
f_b(ax)=f(bax)=f(ab[a,b]x)=h(b[a,b]x)\cdot f(b[a,b]x).
\end{equation}
Since $h$ is a character of $G$ we have
\begin{equation}\label{E:f2}
h(b[a,b]x)=\lambda_1 h(x).
\end{equation}
Using that  $[a,b]$ belongs to the center of  $G$ and
\eqref{E:f_c} we find
\begin{equation}\label{E:f3}
f(b[a,b]x)=f([b,a]bx)=f_{[a,b]}(bx)=\lambda_2 f(bx).
\end{equation}
 Putting together equations  \eqref{E:f1}, \eqref{E:f2}, \eqref{E:f3},
we find
$$
T_af_b=h_1\cdot  f_b
$$
where $h_1=\lambda_1 \lambda_2h\in \mathcal{E}_1(T_a)$. This
completes the proof.
\end{proof}

\begin{proposition}\label{P:nilsystem}
Let  $(X=G/\Gamma, \mathcal{G}/\Gamma,m,T_a)$ be an ergodic
nilsystem and suppose that $X$ is connected. Then $f\in
\mathcal{A}_2(T_a)$ if and only if $f\in L^{\infty}(m)$ and $f$
factors through $G/(G_3 [G_0,G_0]\Gamma)$.
\end{proposition}
\begin{proof}
Suppose first that $f\in L^\infty(m)$ factors through
$G/(G_3[G_0,G_0]\Gamma)$. Replacing $G$ by the group
$G/(G_3[G_0,G_0])$ we can assume that $G$ is $2$-step nilpotent
and that $G_0$ is abelian. In this case, by Theorem~\ref{T:FK}
below
the system is isomorphic to a $2$-step nilpotent affine
transformation on some finite dimensional torus. For such systems
it is easy to verify that $\mathcal{A}_2(T_a)=L^\infty(m)$, and so
$f\in \mathcal{A}_2(T_a)$.

We move now to the converse.   It suffices to show that if  $f\in
\mathcal{E}_2(T_a)$ then $f$ factors through
$G/(G_3[G_0,G_0]\Gamma)$.  We know from \cite{A} that
$f=\text{const}$, so we can assume that $|f|=1$ in which case we
have that $\bar{f}=f^{-1}$. Since $\mathcal{A}_2(T_a)$ is a factor
of $\mathcal{Z}_2(T_a)$, by Theorem~\ref{T:Kronecker} the function
$f$ factors through $G/(G_3\Gamma)$. So it remains to show that
$f$ is invariant under elements in $[G_0,G_0]$. We define the map
$\phi\colon G\to L^\infty(m) $ by $\phi(b)=f(bx)\cdot\bar{f}(x)$.
We need to show that $\phi([G_0,G_0])=1$. First notice that by
Lemma~\ref{L:G_k} the map  $\phi$ takes values in
$\mathcal{E}_2(T_a)$.
Next we claim  that $\phi(G_0)\subset C$ where $C$ is the set of
constant functions of modulus $1$. We will  use a connectedness
argument to show this (similar to the one used in
Lemma~\ref{L:G_k}). If we equip $\mathcal{E}_1(T_a)$ with the
$L^2(m)$ topology then the map $\phi$ is continuous. The connected
component of the function $1$ in $\mathcal{E}_2(T_a)$ is the set
$C$. One can see this by using the fact  that if $f\in
\mathcal{E}_2(T_a)$ is nonconstant then $\int f \ dm=0$ (see
\cite{A}), which implies that $\norm{f-c}_{L^2(m)}=\sqrt{2}$ for
$c\in C$. Since $\phi$ is continuous and $\phi(e)=1$ we have that
$\phi(G_0)\subset C$. Now it is easy to check that $\phi\colon
G_0\to C$ satisfies
$$
\phi(ab)=\phi(a)\cdot\phi(b)=\phi(ba), \qquad
\phi(ab)=\phi(ba)\cdot \phi([a,b]), \quad a,b\in G_0,
$$
which implies (obviously $\phi(x)\neq 0$ for $x\in G$) that
$$
\phi([a,b])=1, \quad a,b\in G_0.
$$
This completes the proof.
\end{proof}

\subsection{Polynomial sequences on nilmanifolds}
If $G$ is a nilpotent Lie group,  $a_1,\ldots, a_k \in G$, and
$p_1, \ldots, p_k $ are integer polynomials
$\mathbb{N}^d\to\mathbb{Z}$, a sequence of the form
$$
 g(\n) =
a_1^{p_1(\n)}a_2^{p_2(\n)}\cdots a_k^{p_k(\n)}
$$
  is called a
\emph{polynomial sequence} in $G$. If the polynomials
$p_1(\n),\ldots,p_k(\n)$ are linear then $g(\n)$ is called a
\emph{linear sequence}.
 The following result of Leibman~(\cite{L0},\cite{L1})
gives information about the orbit closure of polynomial sequences
on nilmanifolds and  helps us handle their uniform distribution
properties\footnote{For the purpose of this article, we will say
that a sequence $\{x(\n)\}_{\n\in\N^d}$ is {\em uniformly
distributed } on the nilmanifold $X$ with the Haar measure $m$, if
for every F{\o}lner sequence $\Phi_N$ in $\N^d$ and continuous
$f\colon X\to \C$, we have $\lim_{N\to\infty} \frac{1}{|\Phi_N|}
\sum_{\n\in\Phi_N}f(x(\n))= \int f\ dm$.} by reducing them to
uniform distribution properties on a certain factor:
\begin{theorem}[{\bf Leibman~\cite{L0},\cite{L1}}]\label{T:L}
Let  $X = G/\gG$ be a nilmanifold and $g(\n)$ be a polynomial
sequence in $G$. Define $Z=G/([G_0,G_0]\Gamma)$ and
$Z'=G/([G,G]\Gamma)$ and let $\pi\colon X\to Z$ and $\pi'
\colon X\to Z'$ be the corresponding  natural projections. Then
for every $x\in X$:

(i) There exist sub-nilmanifolds $Y_i=H x_i$ of $X$, where $H$ is
a closed subgroup of $G$  $($depending on $x$$)$ and
$x_1,\ldots,x_k\in X$, such that
$\overline{\{g(\n)x\}}_{\n\in\mathbb{N}^d}=\bigcup_{i=1}^k Y_i$,
and for $i=1,\ldots, k$ the sequence  $\{g(k\n+i)\}_{n\in \N}$ is
uniformly distributed on $Y_i$. If $g(n)=a^n$ for some $a\in G$,
then $k=1$.

 (ii) If $X$ is connected then  the sequence $\{g(\n)x\}_{\n\in\mathbb{N}^d}$ is
 uniformly distributed in    $X$ if and only
if it is  dense in  $X$ if and only if
$\{g(\n)\pi(x)\}_{\n\in\mathbb{N}^d}$ is dense in $Z$.

(iii) If $X$ is connected and  $a_1,\ldots, a_k\in G$ are
commuting elements that together with $G_0$ generate $G$,   and
$g(\n)=a_1^{p_1(\n)}\cdots a_k^{p_k(\n)}$ is a linear sequence,
then $\{g(\n)x\}_{\n\in\mathbb{N}^d}$ is dense in $X$ if and only
if $\{g(\n)\pi'(x)\}_{\n\in\mathbb{N}^d}$ is dense in $Z'$.
\end{theorem}
We remark that the groups $G_0$, $[G_0,G_0]$, and $[G.G]$ are  normal
subgroups of $G$. Also note that the connected component the identity
element  of the  group $G/[G_0,G_0]$ is abelian. The next result
shows  how we can use this  property
 to our
advantage. In order to state it we need some notation. If $G$ is a
group then a map $T\colon G\to G$ is said to be \emph{affine} if
$T(g) = bA(g)$ for a homomorphism $A$ of $G$ and some $b\in G$.
The homomorphism $A$ is said to be {\em $k$-step nilpotent} if
there exists $k\in\N$ so that $(A-{\text Id})^{k}=0$.
\begin{theorem}[{\bf Frantzikinakis \& Kra~\cite{FK}}]\label{T:FK}
Let $X=G/\gG$ be a connected $k$-step nilmanifold such that  $G_0$ is
abelian. Then every  nilrotation $T_a(x)=ax$ defined on $X$ with
the Haar measure $m$ is isomorphic to a $k$-step nilpotent affine
transformation on some finite dimensional torus. Furthermore, the
conjugation can be chosen to be continuous.
\end{theorem}

Next we give  two applications of Theorems~\ref{T:L} and
\ref{T:FK} that  will be needed in the sequel. We will use the
first one frequently, for example  in the proofs of Theorems~A, B
and F. The simple argument given below was communicated to us by
S. Leibman:
\begin{proposition}\label{L:p(n)}
Suppose that $X=G/\Gamma$ is a nilmanifold, $g\colon \N\to G$ is a
polynomial sequence, and $x\in X$ is such that
$Y=\overline{\{g(n)x\}}_{n\in \N}$ is connected. Then for every
nonconstant integer  polynomial $p$  we have
$Y=\overline{\{g(p(n))x\}}_{n\in \N}$.
\end{proposition}
\begin{proof}
Since $Y$ is connected we have by part $(i)$ of Theorem~\ref{T:L}
 that $Y$ is isomorphic to a subnilmanifold $H/\Delta$ of $X$.
 Hence,
we can assume that $Y=H/\Delta$. By Theorem~\ref{T:L} it suffices
to show that the sequence $\{g(p(n))\pi(x)\}_{n\in \N}$ is
uniformly distributed on $Z=H/([H_0,H_0]\Delta)$ where $\pi \colon
Y\to Z$ is the natural projection. Substituting $H/[H_0,H_0]$ for
$H$ we can assume that $H_0$ is abelian. Suppose that
$g(n)=a_1^{p_1(n)}a_2^{p_2(n)} \cdots a_k^{p_k(n)}$ where
$a_1,\ldots,a_k\in G$. Since $Y$ is connected and $H_0$ is
abelian, by Theorem~\ref{T:FK} we can assume that $Y=\T^m$ and the
nilrotations $T_{a_i}$, $i=1,\ldots,k$, are nilpotent affine
transformations  on $\T^m$. Then the coordinates of the sequence
$\{g(n)x\}_{n\in\N}$ are polynomials in $n$ with real
coefficients, and our problem reduces to the following one: If
$u\colon \N\to \T^m$ is a sequence with polynomial coordinates
such that $\overline{\{u(n)\}}_{n\in\N}=\T^m$, then
$\overline{\{u(p(n))\}}_{n\in\N}=\T^m$ for every nonconstant
polynomial $p$. To see this, first notice that $u$ has the form
$$
u(n)=u_0(n)q+ u_1(n)a_1+\ldots+u_l(n)a_l
$$
where $u_i$ are integer-vector-valued polynomials,
$q\in\mathbb{Q}$, and $a_1,\ldots,a_k$ are linearly independent
irrational numbers. Then  using Corollary 2.4 in \cite{BLL} we
have that $u(n)$ is dense in $\T^m$ if and only if
$$
\text{Span}(u_1(n))+\ldots+\text{Span}(u_k(n)) \! \pmod{1}=\T^m,
$$
where for $u(n)=(q_1(n),\ldots q_r(n))$ we define
$\text{Span}(u(n))=\text{Span}\{(q_1(x),\ldots q_r(x)), x\in\R
\}$. But clearly the last  identity remains  valid if we replace
$n$ with any nonconstant polynomial $p(n)$. This completes the
proof.
\end{proof}
The next lemma   will  be used in the proof of Theorem~E.
\begin{lemma}\label{L:squares}
Suppose that $X=G/\Gamma$ is a nilmanifold, $g\colon \N\to G$ is a
polynomial sequence, and $p$ is an integer polynomial with
$\deg{p}>1$.
Then for every $x\in X$ and $\beta\in\T$ irrational we have
\begin{equation}\label{E:mbm}
\overline{\{(n\beta,g(p(n))x)\}}_{n\in\N}= \mathbb{T}\times
\overline{\{g(p(n)x)\}}_{n\in\N}.
\end{equation}
\end{lemma}
\begin{proof}
Suppose first that the set $Y=\overline{\{g(p(n)x)\}}_{n\in\N}$ is
connected. Working with  the  sequence
$\{h(n)\}_{n\in\N}=\{(n\beta,g(n)x)\}_{n\in\N} $ on $\T\times Y$
and repeating the argument used  in
 the previous lemma we can reduce our problem to the following one:
If $u\colon \N\to \T^m$ is a sequence with polynomial coordinates
such that $\overline{\{u(p(n))\}}_{n\in\N}=\T^m$ and $\deg{p}>1$,
then $\overline{\{(n\beta,u(p(n))\}}_{n\in\N}=\T^{m+1}$. To see
this, first notice that $u$ has the form
$$
u(n)=u_0(n)q+ u_1(n)a_1+\ldots+u_l(n)a_l
$$
where $u_i$ are integer-vector-valued polynomials,
$q\in\mathbb{Q}$, and $a_1,\ldots,a_k$ are rationally independent
 numbers. Since $\overline{\{u(p(n))\}}_{n\in\N}=\T^m$, by
Corollary 2.4 in \cite{BLL} we have
\begin{equation}\label{imu}
\text{Span}\big(u_1(p(n))\big)+\ldots+\text{Span}\big(u_k(p(n))\big)\!\pmod{1}=\T^m.
\end{equation}
We also have
$$
(n\beta,g(p(n))x)=(n,{\bf 0})\beta+(0,u_0(p(n)))q+
(0,u_1(p(n)))a_1+\ldots+ (0,u_l(p(n)))a_l.
$$
Since the polynomials $u_i(p(n))$, $i=1,\ldots,l$, have degree
greater than $1$  it is easy to check that the set
$\overline{\{(n\beta,u(p(n)))\}}_{n\in\N}$ is equal to
$$
 \text{Span}((n,{\bf
0}))+\text{Span}((0,u_1(p(n))))+\ldots+\text{Span}((0,u_k(p(n)))\!
\pmod{1},
$$
and by \eqref{imu} this is equal to $\T^{m+1}$. This completes the
proof.

In the general case we argue as follows: By \cite{L0} there exists
an $r\in\N$ such that $\overline{\{g(p(rn+i))\}}_{n\in\N}$ is
connected for $i=0,\ldots,r-1$. Repeating the previous argument
for the sequence $\{h(rn+i)\}_{n\in\N}$ we find that
$$
\overline{\{h(rn+i)\}}_{n\in\N}= \mathbb{T}\times
\overline{\{g(p(rn+i))x)\}}_{n\in\N}
$$
for $i=0,\ldots,r-1$. This implies \eqref{E:mbm} and completes the
proof.
\end{proof}

\subsection{Limit formula for linear sequences}
In the case where all the polynomials are linear the limit of the
corresponding multiple ergodic averages $(P)$ was computed in
\cite{Z} (for a simpler proof see \cite{BHK}). To state the result
we need some notation. Let $G/\Gamma$ be a  nilmanifold.
Given $l_1,\ldots,l_k\in \N$ define the set
$$
H=\{(g_1^{\binom{l_1}{1}}g_2^{\binom{l_1}{2}} \cdots
g_{k-1}^{\binom{l_1}{k-1}}f_1,\ldots,g_1^{\binom{l_k}{1}}g_2^{\binom{l_k}{2}}
\cdots g_{k-1}^{\binom{l_k}{k-1}}f_k)\colon g_i\in G_i, f_i\in
G_k\}
$$
($\binom{a}{b}=0$ if $a<b$) and let $\Delta=\Gamma^k\cap H$. It
was shown in
 \cite{L}  that $H$ is a closed  subgroup of
$G^k$. The discrete subgroup  $\Delta$ is  cocompact so the
nilmanifold $H/\Delta$ carries a Haar measure, call it $m_{H}$.
The next result is a straightforward generalization of a formula
given by Ziegler in \cite{Z}  that can be obtained using some computations of Leibman in
\cite{L4}:
\begin{theorem}[\bf{Ziegler~\cite{Z}}] \label{T:Z}
Let $(X=G/\gG,\mathcal{G}/\Gamma,T_a,m)$ be an ergodic
nilsystem and $l_1,\ldots,l_k\in \N$. If $f_1,\ldots,f_k\in
L^\infty( X)$ then for a.e. $x=g\Gamma\in X$ we have
\begin{multline*}
\lim_{N-M\to\infty}  \frac{1}{N-M} \sum_{n=M}^{N-1}
 f_1(T_a^{l_1n}x) \cdot \ldots  \cdot  f_k(T_a^{l_kn}x)  =
\int_{H/\Delta} \!   f_1(gy_1\Gamma)\cdot \ldots \cdot
f_k(gy_k\Gamma) \ dm_{H}(y\Delta),
\end{multline*}
where $y=(y_1,\ldots,y_k)$, and $H$, $\Delta$ are as before.
\end{theorem}
Combining this with Theorem~\ref{T:L} we  easily deduce the
following:
\begin{corollary}\label{C:connected}
Let $(X=G/\gG,\mathcal{G}/\Gamma,T_a,m)$ be an ergodic
nilsystem with $X$ connected $($or equivalently $T_a$ is totally
ergodic$)$ and   $l_1,\ldots,l_k\in \Z$. Then for a.e. $x\in X$
the set $H_x=\overline{\{(a^{l_1n}x,a^{l_2n}x,\ldots,
a^{l_kn}x)\}}_{n\in\N}$ is connected.
\end{corollary}
\begin{proof}
By Theorems~\ref{T:L} and \ref{T:Z}   we have that for a.e. $x\in
X$ the set $H_x$ is homeomorphic to the nilmanifold $H/\Delta$
where the subgroup $H$ and $\Delta$ are as before. Since
$X=G/\Gamma$ is connected and $G_i$ is connected for $i\geq 2$ it
follows that $H/\Delta$ is connected. Hence, $H_x $ is connected
for a.e. $x\in X$.
\end{proof}


\section{Weyl complexity for families of three polynomials}\label{S:weyl}
Following \cite{BLL}, we will define the Weyl complexity of a
family $P=\{p_1,\ldots,p_k\}$ of essentially distinct integer
polynomials. Roughly speaking, this notion is designed to capture
the minimum $m\in \N$ for which the factor $\mathcal{Z}_{m-1}$ is
characteristic for the corresponding ergodic averages $(P)$. In
Proposition~\ref{T:weyl} we will give an effective way of
determining the Weyl complexity of any family of three
polynomials.

\subsection{Definition of Weyl complexity and basic
properties}\label{SS:weyl} A \emph{connected Weyl system} is a
system induced by an ergodic nilpotent affine transformation
acting on some finite dimensional torus with the Haar measure. A
\emph{standard Weyl system of level $d$} is a  system induced by a
transformation $T\colon \T^d\to\T^d$ given by
\begin{equation}
\label{E:standard}
T(x_1,\ldots,x_d)=(x_1+\alpha,x_2+x_1,\ldots,x_d+x_{d-1})
\end{equation}
for some irrational number $\alpha \in \T$. A \emph{quasi-standard
Weyl system of level $d$} is a system induced by a transformation
$T\colon \T^d\to\T^d$ given by
\begin{equation}
\label{E:qstandard}
T(x_1,x_2,\ldots,x_d)=(x_1+\alpha_1,x_2+m_{2,1}x_1+\alpha_2,\ldots,x_d+
\sum_{j=1}^{d-1}m_{d,j}x_{j}+\alpha_d)
\end{equation}
where $\alpha_i\in \T$, $\alpha_1$ is irrational, and
$m_{i,i-1}\neq 0$ for all $i=2,\ldots,d$. Note that every
quasi-standard Weyl system is ergodic (\cite{Fu2} page $67$).

 Given a
system $(X,T)$ we denote  the diagonal in $X^{k+1}$ by
$\Delta_{X^{k+1}}$, and we define the orbit of a polynomial family
$P=\{p_1,\ldots,p_k\}$ with respect to the system $(X,T)$ to be
$$
\mathcal{O}(P,\Delta_{X^{k+1}},T)=\{
(x,T^{p_1(n)}x,\ldots,T^{p_k(n)}x)\colon x\in X, n\in \N\}.
$$
\begin{definition}
Let $P=\{p_1,\ldots,p_k\}$ be a family of distinct  integer
polynomials with $p_i(0)=0$ for $i=1,\ldots,k$. The \emph{Weyl
complexity} $W(P)$ is the minimal $r\in \N$ with the following
property: For every $d\in \N$ with $d\geq r$, for
some/every\footnote{It was shown in \cite{BLL} that if
\eqref{E:weeyl} holds for some quasi-standard Weyl system of level
$d$ then it holds for all.} quasi-standard Weyl system $(X,T)$ of
level $d$ we have
\begin{equation}\label{E:weeyl}
\{(0,\ldots,0,x_r,\ldots,x_d)\}^{k+1}\subset
\overline{\mathcal{O}(P,\Delta_{X^{k+1}},T)}.
\end{equation}
For a general family $P=\{p_1,\ldots,p_k\}$ of essentially
distinct
 polynomials  we define $W(P)=W(p_1-p_1(0),\ldots,p_k-p_k(0))$.
\end{definition}

The next two results give  equivalent characterizations of the
Weyl complexity that are better suited for our purposes. The first
follows easily from the definition.
\begin{proposition}\label{P:weyl1}
The Weyl complexity $W(P)$ of a  family $P=\{p_1,\ldots,p_k\}$ of
essentially distinct integer polynomials is the maximal $s\in\N$
$($or $1$ if there is no such $s$$)$ with the following property:
For some/every quasi-standard Weyl system $(X,T)$ of level $s-1$
of the form \eqref{E:qstandard}, there exist characters $\chi_i$
of $X$, $i=0,\ldots,k$, at least one of which depends nontrivially
on the variable $x_{s-1}$, such that
$$
\chi_0(x)\cdot \chi_1(T^{p_1(n)}x)\cdot \ldots\cdot
\chi_k(T^{p_k(n)}x)=1
$$
for every $x\in \T^{s-1}$.
\end{proposition}
For a proof of the next result see the remarks after Proposition
$5.1$ in \cite{BLL}.
\begin{proposition}\label{P:weyl3}
The Weyl complexity $W(P)$ of a  family of essentially distinct
integer polynomials $P=\{p_1,\ldots,p_k\}$  is the minimal $m\in
\N$ with the following property: For every connected Weyl system
$(X,T)$ the factor $\mathcal{Z}_{m-1}$ is characteristic for
$L^2$-convergence or weak convergence of the averages $(P)$.
\end{proposition}

We remark that for a quasi-standard Weyl system  of the form
\eqref{E:qstandard}  the factor $\mathcal{Z}_m$ coincides with the
sub-$\sigma$-algebra of sets that depend only on  the first $m$
coordinates.

We will make frequent use of the following simple identity:
\begin{proposition}\label{P:weyl4}
If $\{p_1,\ldots,p_k\}$ is a family of essentially distinct
polynomials then
$$
W(p_1,p_2,\ldots,p_k)=W(p_1-p_k,p_2-p_k,\ldots,p_{k-1}-p_k,-p_k).
$$
\end{proposition}
\begin{proof}
 This follows immediately from Proposition~\ref{P:weyl3} and the
identity
$$
\int f_0 \cdot T^{p_1(n)}f_1\cdot\ldots\cdot T^{p_k(n)}f_k \ d\mu=
\int T^{-p_k(n)}f_0 \cdot T^{p_1(n)-p_k(n)}f_1\cdot\ldots\cdot
T^{p_{k-1}(n)-p_k(n)}f_{k-1}\cdot f_k \ d\mu.
$$
\end{proof}

\subsection{Different scenarios for the Weyl complexity of three  polynomials}
We will give an explicit criterion for determining
$W(p_1,p_2,p_3)$. We first show:
\begin{proposition}\label{<=3}
If $\{p_1,p_2,p_3\}$ is a family of essentially distinct
polynomials then  $W(p_1,p_2,p_3)\leq 3$.
\end{proposition}
\begin{proof}
 We argue by contradiction. Suppose that $W(p_1,p_2,p_3)\geq 4$.
We can assume that $p_i(0)=0$ for $i=1,2,3$. Consider the
quasi-standard Weyl system $(\T^3,T)$ where
$$
T(x_1,x_2,x_3)=(x_1+\alpha,x_2+2x_1+\alpha,x_3+3x_1+3x_2+\alpha).
$$
By Proposition~\ref{P:weyl1} there exist characters
$\chi_0,\chi_1,\chi_2,\chi_3$ of $\T^3$,  at least one of which
depends nontrivially on the variable $x_3$, such that
\begin{equation}\label{E:766}
 \chi_0(x)\cdot
\chi_1(T^{p_1(n)}x)\cdot \chi_2(T^{p_2(n)}x) \cdot
\chi_3(T^{p_3(n)}x)=1
\end{equation}
for all $x\in \T^3$. We use that
$$
T^n(x_1,x_2,x_3)=(x_1+n\alpha,x_2+2nx_1+n^2\alpha,x_3+3nx_2+3n^2x_1+n^3\alpha)
$$
and substitute in  \eqref{E:766}. Suppose that
$$
\chi_i(x_1,x_2,x_3)=e(k_ix_1+l_ix_2+m_ix_3)
$$
for some integers $k_i,l_i,m_i$ for  $i=0,1,2,3$.  Plugging in
\eqref{E:766} we get that the system
\begin{align}
\label{E:1}&k_1p_1+l_1p_2+m_1p_3&=0\\
\label{E:2}&3k_1p_1^2+3l_1p_2^2+3m_1p_3^2+2k_2p_1+2l_2p_2+2m_2p_3&=0\\
\label{E:3}&k_1p_1^3+l_1p_2^3+m_1p_3^3+k_2p_1^2+l_2p_2^2+m_2p_3^2+
k_3p_1+l_3p_2+m_3p_3& =0
\end{align}
has a  solution on  the integers $k_i,l_i,m_i$, $i=1,2,3$, with at
least one of the $k_1,l_1,m_1$   nonzero. Let $d_i=\deg{p_i}$,
$i=1,2,3$, and $a_1,b_1,c_1$ be the leading coefficients of the
polynomials $p_1,p_2,p_3$. After rearranging the polynomials we
can assume that $d_1\geq d_2\geq d_3$. We consider three cases:

{\bf Case 1.} If $d_1>d_2\geq d_3$ then \eqref{E:1} gives $k_1=0$,
so $d_2=d_3=d>0$. If $k_2=0$ then looking at the leading
coefficients of the polynomials in  \eqref{E:2} we  get that
$l_1=-m_1$ which implies (using \eqref{E:1}) that $p_2=p_3$, a
contradiction. Hence, $k_2\neq 0$ and  since $l_1\neq -m_1$ we get
from \eqref{E:2} that $d_1=2d$. But then the polynomial on the
left hand side of \eqref{E:3} has degree $4d$, a contradiction.

{\bf Case 2.} If $d_1\geq d_2> d_3$ then using equations
\eqref{E:1}, \eqref{E:2} (which are the same as \eqref{E:(1)},
\eqref{E:(2)}), and arguing as in Case $2$ of
Proposition~\ref{T:weyl} below, we derive that $p_1=p_2$, a
contradiction.

{\bf Case 3.} If $d_1=d_2=d_3=d$ then looking at the leading
coefficients of the polynomials in \eqref{E:1}, \eqref{E:2},
\eqref{E:3}, we get that the system
\begin{align*}
k_1a_1+l_1b_1+m_1c_1&=0\\
k_1a^2_1+l_1b^2_1+m_1c^2_1&=0\\
 k_1a^3_1+l_1b^3_1+m_1c^3_1&=0
\end{align*}
has a nontrivial integer solution on $k_1,l_1,m_1$. The
determinant of the corresponding matrix  is
$a_1b_1c_1(a_1-b_1)(b_1-c_1)(c_1-a_1)$. Since $a_1,b_1,c_1$ are
nonzero, two of them must be  equal. Without loss of generality we
can assume that  $a_1=b_1$.  Then after replacing $p_1$ with
$q_1=-p_1$, $p_2$ with $q_2=p_2-p_1$, $p_3$ with $q_3=p_3-p_1$,
and using Proposition~\ref{P:weyl4}, reduces our problem to either
Case $1$ or Case $2$. So we get again a contradiction showing that
$W(p_1,p_2,p_3)\leq 3$.
\end{proof}
We will also need the following simple lemma:
\begin{lemma}\label{L:system}
Suppose that  $a_1,b_1,c_1\in\mathbb{Z}$ are nonzero and distinct,
$a_2,b_2,c_2\in\mathbb{Z}$, and
\begin{align*}
k_1a_1+l_1b_1+m_1c_1&=0, \qquad \qquad \qquad\ \ \!   k_1a_2+l_1b_2+m_1c_2=0\\
k_1a^2_1+l_1b^2_1+m_1c^2_1&=0,\qquad \qquad
k_1a_1a_2+l_1b_1b_2+m_1c_1c_2=0,
\end{align*}
for some integers $k_1,l_1,m_1$, not all of them zero. Then there
exist $r,s\in\mathbb{Q}$ such that
$$
(a_1,a_2)=r(b_1,b_2)=s(c_1,c_2).
$$
\end{lemma}
\begin{proof}
Without loss of generality we can assume that $l_1\neq 0$.
Performing some elementary operations we get the system
\begin{align*}
(a_1b_1-b_1^2)l_1+(a_1c_1-c_1^2)m_1&=0\\
(b_1a_2-b_1b_2)l_1+(c_1a_2-c_1c_2)m_1&=0\\
(b_2a_1-b_1b_2)l_1+(c_2a_1-c_1c_2)m_1&=0.
\end{align*}
Using that  $b_1c_1\neq 0$ the first two equations easily imply
that   $a_1c_2=c_1a_2$. Since $(a_1- b_1)(a_1-c_1)\neq 0$ the
first and third equation easily imply  that $b_1c_2=b_2c_1$. The
result follows.
\end{proof}

We can now prove the main result of this section:
\begin{proposition}\label{T:weyl}
Let $p_1,p_2,p_3$ be essentially distinct polynomials and let
$\tilde{p}_i=p_i-p_i(0)$ for  $i=1,2,3$. Then:

\noindent $(i)$ $W(p_1,p_2,p_3)=1$ if and only if the polynomials
$\tilde{p}_1,\tilde{p}_2,\tilde{p}_3$ are linearly independent.

\noindent $(ii)$ $W(p_1,p_2,p_3)=3$ if and only if  some
permutation of the polynomials
$\tilde{p}_1,\tilde{p}_2,\tilde{p}_3$ has the  form
$$ (a) \  (lp,mp,kp^2+rp),\quad \text{or}\quad
(b) \  (kp^2+lp,kp^2+mp,kp^2+rp),
$$
for some  $k,l,m,r\in\Z$.

\noindent $(iii)$ In all other cases  $W(p_1,p_2,p_3)=2$.
\end{proposition}
\begin{proof}
We can assume that $p_i(0)=0$ for $i=1,2,3$. We first show part
$(i)$. Consider the standard Weyl system $(\T,T)$ of level $1$
induced by the transformation $Tx=x+\alpha$, where $\alpha\in \R$
is irrational. Let $\chi_i(x)=e(m_ix)$, where $m_i\in \Z$ for
$i=0,1,2,3$, be characters of $\T$. Since the polynomials
$p_1,p_2, p_3$ are linearly independent the equation
$$
\chi_0(x)\cdot \chi_1(T^{p_1(n)}x)\cdot \chi_2(T^{p_2(n)}x)\cdot
\chi_3(T^{p_3(n)}x)=1
$$
gives that $m_i=0$ for $i=0,1,2,3$. By Proposition~\ref{P:weyl1}
we have  that $W(p_1,p_2,p_3)=1$.

To show part $(i)$ we first notice that since by
Proposition~\ref{<=3} we have $W(p_1,p_2,p_3)\leq 3$, it remains
to show that $W(p_1,p_2,p_3)\geq 3$ if and only if the polynomials
have the form $(a)$ or $(b)$. To do this consider the
quasi-standard Weyl system $(\T^2,T)$ defined by
$$
T(x_1,x_2)=(x_1+\alpha,x_2+2x_1+\alpha).
$$
By Proposition~\ref{P:weyl1} we have $W(p_1,p_2,p_3)\geq 3$ if and
only if there exist characters $\chi_0,\chi_1,\chi_2,\chi_3$ of
$\T^2$, at least one of which depends nontrivially on the variable
$x_2$, such that
\begin{equation}\label{E:76}
 \chi_0(x)\cdot
\chi_1(T^{p_1(n)}x)\cdot \chi_2(T^{p_2(n)}x) \cdot
\chi_3(T^{p_3(n)}x)=1
\end{equation}
for all $x\in \T^3$. We use that
$$
T^n(x_1,x_2)=(x_1+n\alpha,x_2+2nx_1+n^2\alpha)
$$
and substitute in  \eqref{E:76}. We get  that $W(p_1,p_2,p_3)\geq
3$ if and only if  the system
\begin{align}
\label{E:(1)}k_1p_1+l_1p_2+m_1p_3 \qquad \qquad  \qquad \qquad
\quad
&=0\\
\label{E:(2)}k_1p_1^2+l_1p_2^2+m_1p_3^2+k_2p_1+l_2p_2+m_2p_3&=0
\end{align}
has an integer solution on the $k_i,l_i,m_i$, $i=1,2$, with at
least one of the $k_1,l_1,m_1$ nonzero.

If the polynomial family  has the form $(a)$ then the following
are  eligible  solutions to the previous system: $(i)$ If $k\neq
0$ then  $k_1=mk$, $l_1=-kl$, $m_1=0$, $k_2=m(l-m)r$, $l_2=0$,
$m_2=ml(l-m)$, $(ii)$ If  $k=0$ then $k_1=rm(l-r)$, $k_2=rl(r-m)$,
$k_3=lm(m-l)$, $k_2=0$, $l_2=-r$, $m_2=m$. Hence,
$W(p_1,p_2,p_3)=3$. By Proposition~\ref{P:weyl4} we get that the
same is true for any polynomial family of the form $(b)$.

We now  focus on the hardest part of the result which is to show
that if  $W(p_1,p_2,p_3)= 3$ then some permutation of the
polynomials has either the form $(a)$ or $(b)$. Let
$$
p_1(n)=a_1n^{d_1}+\ldots+a_{d_1}n,\quad
p_2(n)=b_1n^{d_2}+\ldots+b_{d_2}n,\quad
p_3(n)=c_1n^{d_3}+\ldots+c_{d_3}n
 $$
for some $d_i\in\mathbb{N}$, $i=1,2,3$, and
$a_i,b_i,c_i\in\mathbb{Z}$ with $a_1,b_1,c_1\neq 0$.
  After rearranging the
polynomials we can assume that $d_1\geq d_2\geq d_3$.
We consider the following three cases:

{\bf Case 1.} If  $d_1>d_2\geq d_3$ we will show that some
permutation of the polynomials has either the form  $(a)$ or
$(b)$. From \eqref{E:(1)} we get $k_1=0$, so $p_2,p_3$ are integer
multiples of the same integer polynomial $p$. Using this and
\eqref{E:(2)}, we get that $k_2p_1$ is an integer combination of
$p$ and $p^2$. This easily implies that the polynomials have the
form  $(a)$, possibly with some rational multiple of $p$ in place
of $p$.

{\bf Case 2.} If $d_1\geq d_2> d_3$ we will show that $p_1=p_2$, a
contradiction. From  \eqref{E:(1)} we get that $d_1=d_2=d$ and
looking at the leading coefficients of the polynomials in
$\eqref{E:(1)}$ and $\eqref{E:(2)}$ we get the system
$$
k_1a_1+l_1b_1=0, \quad k_1 a^2_1+l_1b^2_1=0.
$$
Since $a_1,b_1\neq 0$ and  $k_1,k_2$ are not both zero, we easily
get that $a_1=b_1$ and   $k_1+l_1=0$. If  $d_0=\max\{j\in
\mathbb{N}\colon a_i=b_i, \text{ for } 1\leq i\leq j\}$ (it is
well defined since $a_1=b_1$), it suffices to show that $d_0=d$.
Suppose not, then we can write $p_1=q+p_1'$, $p_2=q+p_2'$, where
$q(n)=a_1n^d+\ldots+a_{d_0}n^{d_0-k+1}$, the degrees $d'_1,d'_2$
of $p'_1,p_2'$ are not both  zero, and they do not exceed $d_0-k$.
By possibly permuting the polynomials we can further assume that
\begin{equation}\label{E:degrees}
 d'_2\leq d'_1<d, \quad \deg{(p_1'-p_2')}=d'_1.
\end{equation}
Substituting $p_1=q+p_1'$, $p_2=q+p_2'$, and $k_1=-l_1$ in
equations \eqref{E:(1)} and \eqref{E:(2)} gives the system
\begin{align}
\label{E:(1')} k_1p_1'-l_1 p_2'+m_1p_3 \qquad \qquad  \qquad
\qquad \qquad \qquad  \qquad \qquad \quad
&=0\\
\label{E:(2')}
k_1p_1'^2-k_1p_2'^2+m_1p_3^2+2k_1q(p_1'-p_2')+k_2p_1+l_2p_2+m_2p_3
&=0.
\end{align}
By \eqref{E:(1')} we get $m_1\neq 0$ (otherwise $k_1=m_1=0$)  and
the polynomial $p_3$ has degree at most $d_1'$.  By
\eqref{E:degrees} the polynomial $q(p_1'-p_2')$ has degree
$d+d_1'$ which is greater than the degree of all other polynomials
that appear in \eqref{E:(2')}. This can only happen if $k_1=0$,
which gives $m_1=0$, contradicting our assumption that one of the
integers $k_1,l_1,m_1$ is nonzero. Hence,
 $d_0=d$ which gives that
$p_1=p_2$.

{\bf Case 3.} If $d_1=d_2=d_3=d$ we will show that the polynomials
have the form $(b)$. We consider two subcases. Suppose that two of
the three leading coefficients are the same, say for example that
$a_1=b_1$ (the other cases can be treated similarly). Then after
replacing $p_1$ with $q_1=-p_1$, $p_2$ with $q_2=p_2-p_1$, $p_3$
with $q_3=p_3-p_1$, and using Proposition~\ref{P:weyl4}, reduces
our problem to either Case $1$ or Case $2$. Since Case 2 is
impossible, the polynomials $q_1,q_2,q_3$ have the form $(a)$. It
follows that the polynomials $p_1,p_2,p_3$ have the form $(b)$ for
some $k\neq 0$.

 So it remains to deal with the case where all three polynomials have degree $d$ and
 their leading coefficients $a_1,b_1,c_1$ are distinct.  In this
 case we will show
that the polynomials have the form $(b)$ with $k=0$. The case
where  $d=1$ is trivial so we can assume that $d\geq 2$. There
exist nonzero $r,s\in\mathbb{Q}$ such that $ a_1=rb_1=sc_1. $ We
will show by induction on $t$ that for all $1\leq t\leq d$ we have
\begin{equation}\label{E:vector}
(a_1,a_2,\ldots,a_t)=r(b_1,b_2,\ldots,b_t)=s(c_1,c_2,\ldots,c_t).
\end{equation}
The $t=d$ case gives that the polynomials $p_1,p_2,p_3$ have the
form $(b)$ with $k=0$. For $t=1$ the statement is true by
assumption. To better illustrate the idea of the inductive step we
first work out the $t=2$ case. Looking at the coefficient of $n^d$
and $n^{d-1}$ in \eqref{E:(1)}, and the coefficient of $n^{2d}$
and $n^{2d-1}$ in \eqref{E:(2)}, we get (for $d\geq 2$ we have
$2d-1>d$) the system
\begin{align*}
k_1a_1+l_1b_1+m_1c_1&=0,   &k_1a_2+l_1b_2+m_1c_2=0\\
k_1a^2_1+l_1b^2_1+m_1c^2_1&=0, &k_1a_1a_2+l_1b_1b_2+m_1c_1c_2=0.
\end{align*}
Since the integers $a_1,b_1,c_1$ are nonzero and distinct we get
by Lemma~\ref{L:system} that $ (a_1,a_2)=r(b_1,b_2)=s(c_1,c_2) $
for some nonzero   $r,s\in \mathbb{Q}$, proving that
\eqref{E:vector} holds for $t=2$.

Inductive step: Suppose that \eqref{E:vector} holds for some
$t\in\N$ with  $1\leq t< d$, we will show that it holds for $t+1$.
So we need to establish  that
\begin{equation}\label{E:ind}
(a_1,a_{t+1})=r(b_1,b_{t+1})=s(c_1,c_{t+1}).
\end{equation}
 Looking at the coefficient
of $n^d$ and $n^{d-t}$ in \eqref{E:(1)}, and the coefficient of
$n^{2d}$  and $n^{2d-t}$ in \eqref{E:(2)}, we get the system
\begin{align}\label{E:system2}
k_1a_1+l_1b_1+m_1c_1&=0,   &k_1a_{t+1}+l_1b_{t+1}+m_1c_{t+1}=0\\
\label{hg4} k_1a^2_1+l_1b^2_1+m_1c^2_1&=0, & \sum_{i=1}^t (
k_1a_ia_{t+2-i}+l_1 b_ib_{t+2-i}+m_1c_ic_{t+2-i})=0.
\end{align}
Since $a_1=rb_1=sc_1$, the first equation in $\eqref{hg4}$ gives
that
\begin{equation}\label{E:rs}
k_1+l_1r^{-2}+m_1s^{-2}=0.
\end{equation}
If  $1\leq i\leq t$,  by the inductive hypothesis  we have
$a_i=rb_i=sc_i$. So for $2\leq i\leq t$ (then $t+2-i\leq t$) we
get
$$
k_1a_ia_{t+2-i}+l_1
b_ib_{t+2-i}+m_1c_ic_{t+2-i}=a_ia_{t+2-i}(k_1+l_1r^{-2}+m_1s^{-2})=0,
$$
where the last equality holds from \eqref{E:rs}. This shows that
in the second equation in \eqref{hg4} all the terms in the sum
except the first one are zero, hence
$$
k_1a_1a_{t+1}+l_1b_1b_{t+1}+m_1c_1c_{t+1}=0.
$$
If we replace the second equation in \eqref{hg4} with this simpler
one,    Lemma~\ref{L:system} applies and gives  \eqref{E:ind}.
This completes the induction and the proof.
\end{proof}

\section{Characteristic factors for  the families $\{l_1p,l_2p,\ldots,l_kp\}$ and  $\{p_1,p_2,p_3\}$}
\label{S:cf}
In this section we will prove Theorems~A and B.

\subsection{Some preparatory work} We start with some preliminary
results. We say that a collection $\mathcal{P}$ of integer
polynomial families is \emph{eligible} if whenever
$\{p_1(n),\ldots,p_k(n)\}\in \mathcal{P}$ then $(i)$
$\{p_1(rn+s),\ldots, p_k(rn+s)\}\in \mathcal{P}$ for every
$r\in\N, s=0,\ldots,r-1$, and $(ii)$
$\{cp_1(n),\ldots,cp_k(n)\}\in \mathcal{P}$ for every nonzero
$c\in\mathbb{Q}$, as long as $cp_i\in\Z[t]$ for $i=1,\ldots,k$.
\begin{proposition}\label{P:te}
Let $\mathcal{P}$ be an eligible collection of $k$-term polynomial
families. Suppose that there exists an $m\in\N$ such that for
every totally ergodic nilsystem and every $\{p_1,\ldots,p_k\}\in
\mathcal{P}$ the factor $\mathcal{Z}_m$ is characteristic for weak
convergence of the ergodic averages $(P)$. Then the same is true
for $L^2$-convergence and for every ergodic system.
\end{proposition}
\begin{proof}
We can assume that $p_i(0)=0$ for $i=1,\ldots,k$. By
Theorem~\ref{T:HKL} we know that the averages $(P)$ converge in
$L^2(\mu)$, so the corresponding weak and strong limits coincide.

Suppose that the factor $\mathcal{Z}_m$ satisfies the assumption
of the Proposition. It suffices to show that for every ergodic
system $(X,\mathcal{X},\mu,T)$ if $f_i\in L^\infty(\mu)$ for
$i=1,\ldots,k$ and  $\E(f_i|\mathcal{Z}_m)=0$ for some
$i=1,\ldots,k$, then the averages $(P)$ converge to $0$ in
$L^2(\mu)$ as $N-M\to\infty$. Without loss of generality we can
assume that $i=1$.   For ergodic systems, by Theorem~\ref{T:HKL}
there exists a characteristic factor that is an inverse limit of
nilsystems induced by some $T$-invariant sub-$\sigma$-algebras
$\{\mathcal{X}_i\}_{j\in\N}$. Since
$\E(f_1|\mathcal{Z}_m(\mathcal{X}))=0$ implies that
$\E(f_1|\mathcal{Z}_m(\mathcal{X}_j))=0$, for $j\in\Z$, an
approximation argument allows us to assume that our system is an
ergodic nilsystem, say  $(X=G/\Gamma,\mathcal{G}/\Gamma,m,T_a)$.
By Proposition~\ref{basic} there exists an $r\in\N$ such that the
ergodic components of $T_a^r$ are totally ergodic. Since
$p_i(0)=0$, we have that $p_i(nr)=rq_i(n)$ for some integer
polynomials $q_i$, for $i=1,\ldots,k$. Because $\mathcal{P}$ is
eligible we have that $\{q_1,\ldots,q_k\}\in \mathcal{P}$. We know
from \cite{L3} that for every nonzero integer $r$ and $m\in\N$ we
have $\cZ_m(T_a)=\cZ_m(T_a^r)$. Since $T^r$ has finitely many
ergodic components, it follows that if $\E(f|\mathcal{Z}_m)=0$
then the same holds for the ergodic components of $T_a^r$. So
using our assumption for the ergodic components of $T_a^r$ and the
polynomial family $\{q_1,\ldots,q_k\}$,  we get that the averages
$(P)$ converge to $0$ in $L^2(\mu)$ as $N-M\to\infty$ if we
substitute $p_i(rn)$ for $p_i(n)$ for $i=1,\ldots,k$. Finally,
since $\E(f_1|\mathcal{Z}_m)=0$ implies that
$\E(T_a^jf_1|\mathcal{Z}_m)=0$, for $j\in\N$, a similar argument
shows that the limit is also zero if we substitute $p_i(nr+s)$ for
$p_i(nr)$ in $(P)$ for $s=0,\ldots,r-1$. It follows that the
averages $(P)$ converge to $0$ in $L^2(\mu)$ as $N-M\to\infty$,
completing the proof.
\end{proof}
Next we prove Theorem~A.
\begin{proof}[Proof of Theorem~A]
Let $p$ be a nonconstant integer polynomial. We  first claim that
for totally ergodic systems the $L^2$-limit of the ergodic
averages associated to the families
$\{l_1p(n),l_2p(n),\ldots,l_kp(n)\}$ and
$\{l_1n,l_2n,\ldots,l_kn\}$ are the same  (a formula for the limit
is then follows from Theorem~\ref{T:Z}). Using Theorem~\ref{T:HKL}
and an approximation argument it suffices to check this for every
totally ergodic nilsystem. So let $(X=G/\Gamma,\mathcal{G}/\Gamma,
m, T_a)$ be such a system. It suffices to show that  for a.e.
$x\in X$ the sequences
$\{(a^{l_1n}x,a^{l_2n}x,\ldots,a^{l_kn}x)\}_{n\in\N}$ and
$\{(a^{l_1p(n)}x,a^{l_2p(n)}x,\ldots,a^{l_kp(n)}x)\}_{n\in\N}$ are
equidistributed, or equivalently that  the sequences
$\{g(n)\tilde{x}\}_{n\in\N}$ and
 $\{g(p(n))\tilde{x}\}_{n\in\N}$ are equidistributed, where
 $g(n)$$=(a^{l_1n},$$a^{l_2n},$$\ldots,$ $ a^{l_kn})$ is a linear sequence in $G^k$ and
 $\tilde{x}=(x,\ldots,x)\in X^k$.
 By Theorem~\ref{T:L} it is enough  to show that
 for a.e. $x\in X$ the two sequences have the same
 closure.   By
 Corollary~\ref{C:connected} the set
 $\overline{\{g(n)x\}}_{n\in\N}$ is connected for a.e.
 $x\in X$, so  Proposition~\ref{L:p(n)} applies and gives the required identity.

We know from  \cite{HK}, \cite{L3} that the factor
$\mathcal{Z}_{k-1}$ is characteristic for the family
$\{l_1n$$,l_2n,$$\ldots,$ $l_kn\}$, hence $\mathcal{Z}_{k-1}$ is
also characteristic for the family $
\{l_1p(n),l_2p(n),\ldots,l_kp(n)\}$ for totally ergodic systems.
Since the collection of polynomial families of the form
$\{l_1p(n),l_2p(n),\ldots,$ $l_kp(n)\}$ with $p\in\Z[t]$
nonconstant is eligible, by Proposition~\ref{P:te} the factor
$\cZ_{k-1}$ is also characteristic for every ergodic system.
 It was shown in \cite{Z1} that  $\mathcal{Z}_{k-1}$ is in fact
 the smallest characteristic factor for the family $\{l_1n,l_2n,\ldots,l_kn\}$,
 the same argument shows that this is also the case  for any family of the
 form $\{l_1p(n),l_2p(n),\ldots,l_kp(n)\}$ where $p\in \Z[t]$ is  nonconstant.
\end{proof}

The next two lemmas will enable us to show that the Kronecker
factor is characteristic for the averages $(P)$ when $k=3$ and
$W(p_1,p_2,p_3)=2$.

\begin{lemma}\label{L:w2'}
Let $k_1,k_2,l_1,l_2\in \Z$ be such that the polynomials
$k_1m,k_2n, l_1m+l_2n$ are   distinct. Then for ergodic systems
the Kronecker factor is characteristic  for $L^2$-convergence of
the averages
\begin{equation}\label{E:double} \frac{1}{(N-M)^2} \sum_{m,n=M}^{N-1}
T^{k_1m}f_1\cdot T^{k_2n}f_2\cdot T^{l_1m+l_2n}f_3.
\end{equation}
\end{lemma}
\begin{proof}
Let $(X,\mathcal{X},\mu,T)$ be an ergodic system and suppose that
$f_1,f_2,f_3\in L^{\infty}(\mu)$ with $\norm{f_i}_\infty\leq 1$
for $i=1,2,3$. It suffices to show that if $\E(f_i|\mathcal{K})=0$
for  some $i=1,2,3$ then the $L^2$-limit of the averages in
\eqref{E:double} is zero. Suppose that $E(f_3|\mathcal{K})=0$, the
argument is identical if  $\E(f_2|\mathcal{K})=0$, and if
$\E(f_1|\mathcal{K})=0$ we only have to interchange the role of
$m$ and $n$.  By Theorem~\ref{T:L} \footnote{More accurately, we
have to combine  Theorem~\ref{T:L} with a result in  \cite{L3}
that reduces the study of the limiting behavior of linear multiple
ergodic averages along any F{\o}lner sequence to nilsystems.} the
$L^2$-limit
\begin{equation}\label{E:double'}
\lim_{N\to\infty} \frac{1}{|F_N|} \sum_{m,n\in F_N}
T^{k_1m}f_1\cdot T^{k_2n}f_2\cdot T^{l_1m+l_2n}f_3
\end{equation}
does not depend on the choice of the F{\o}lner sequence $F_N$. We
claim that the $L^2$-limit  \eqref{E:double'} is zero for the
F{\o}lner sequence
\begin{equation}\label{E:Folner}
F_N=\{0\leq m\leq N,0\leq n\leq a(N)\},
\end{equation}
 where $a(N)$ is an  increasing  sequence of integers that will be chosen later.
 We start by
using the well known fact that for the family $\{an,bn\}$, where
$a,b$ are distinct integers, the Kronecker factor is
characteristic (this is implicit in \cite{Fu1}). Since
$\E(f_3|\mathcal{K})=0$, we get  that for every   $N\in\N$ there
exists an $a(N)\in \N$ such that
\begin{equation}\label{E:a1}
\norm{\frac{1}{a(N)}\sum_{n=0}^{a(N)-1} T^{k_2n}f_2\cdot
T^{l_2n}(T^{l_1m}f_3)}_{L^2(\mu)}\leq \frac{1}{N}
\end{equation}
for all $ m\in \{0,1,\ldots,N-1\}$. Furthermore, we can make sure
that the sequence $a(N)$ is increasing in $N$. We  have
\begin{multline}\label{E:a2}
\norm{ \frac{1}{|F_N|} \sum_{m,n\in F_N} T^{k_1m}f_1\cdot
T^{k_2n}f_2\cdot T^{l_1m+l_2n}f_3}_{L^2(\mu)}\leq\\
\frac{1}{N}\sum_{m=0}^N \norm{\frac{1}{a(N)}\sum_{n=0}^{a(N)-1}
T^{k_2n}f_2\cdot T^{l_2n}(T^{l_1m}f_3)}_{L^2(\mu)}.
\end{multline}
Combining \eqref{E:a1} and \eqref{E:a2}, we get that for the
choice of F{\o}lner sequence made in \eqref{E:Folner} the
$L^2$-limit in \eqref{E:double'} is zero, and so the same is true
for the $L^2$-limit of the averages \eqref{E:double}.
\end{proof}

\begin{lemma}\label{L:w2}
Let $(X,\mathcal{X},\mu,T)$ be a totally ergodic system, $p_1,p_2$
be linearly independent integer polynomials, and
 $k_1$,$k_2$,$l_1$,$l_2$$\in \Z$ be
 such that the  family
 $P=\{k_1p_1,k_2p_2,l_1p_1+l_2p_2\}$ has Weyl complexity $2$.
If $f_0,f_1,f_2,f_3\in L^\infty(\mu)$ then the averages
\begin{equation}\label{b1}
\frac{1}{N-M} \sum_{n=M}^{N-1} \int f_0 \cdot
T^{k_1p_1(n)}f_1\cdot T^{k_2p_2(n)}f_2\cdot
T^{l_1p_1(n)+l_2p_2(n)}f_3\ d\mu
\end{equation}
and
\begin{equation}\label{b2}
\frac{1}{(N-M)^2} \sum_{n,r=M}^{N-1}\int f_0 \cdot
T^{k_1n}f_1\cdot T^{k_2r}f_2\cdot T^{l_1n+l_2r}f_3\ d\mu
\end{equation}
have the same limit as $N-M\to\infty$.
\end{lemma}
\begin{proof}
By Theorem~\ref{T:HKL} and Lemma~\ref{L:w2'} there exists a factor
of the system  that is characteristic  for both averages and is an
inverse limit of finite step nilsystems. So using an approximation
argument it suffices to verify the lemma when the system is a
totally ergodic nilsystem, say
$(X=G/\Gamma,\mathcal{G}/\Gamma,m,T_a)$. By Proposition~\ref{P:te}
the set $X$ is connected so using Theorem~\ref{T:L} it suffices to
show that for every $x\in X$ the sequences
\begin{equation}\label{E:seq1}
\{(a^mx,a^{m+k_1n}x,a^{m+k_2r}x,a^{m+l_1n+l_2r}x)\}_{m,n,r\in\N}
\end{equation}
and
\begin{equation}\label{E:seq2}
\{(a^mx,a^{m+k_1p_1(n)}x,a^{m+k_2p_2(n)}x,a^{m+l_1p_1(n)+l_2p_2(n)}x)\}_{m,n\in\N}
\end{equation}
have the same closure.

Consider the subgroup of $G^4$ defined by
$$
H=\{(g f_1,gh_1^{k_1}f_2,gh_2^{k_2}f_3,gh_1^{l_1}h_2^{l_2}f_4)\in
G^4\colon g,h_1,h_2\in G, f_1,f_2,f_3,f_4\in G_2\},
$$
and  $\Delta=H\cap \Gamma^4$. Using part $(iii)$ of
Theorem~\ref{T:L}  we can show  that the closure of the sequence
in \eqref{E:seq1} is the connected nilmanifold $H/\Delta$, where
$\Delta=H\cap \Gamma^4$ (alternatively we can directly quote a
more general result proved in Section $4$ of \cite{L4}). It
remains to show that the closure of the sequence in \eqref{E:seq2}
is equal to $H/\Delta$ as well. To do this we are going to apply
Theorem~\ref{T:L}. First notice that if $a_0=(a,a,a,a)$,
$a_1=(e,a^{k_1},e,a^{l_1})$, $a_2=(e,e,a^{k_2},a^{l_2})$,
  and $\tilde{x}=(x,x,x,x)$, then the  sequences in \eqref{E:seq1}
  and \eqref{E:seq2} take the form
 $ \{h_1(m,n,r)\tilde{x}\}_{m,n,r\in\N}$ and
  $ \{h_2(m,n)\tilde{x}\}_{m,n\in\N}$,
  where
 $ h_1(m,n,r)=a_0^ma_1^n a_2^r$ and  $h_2(m,n)=a_0^m a_1^{p_1(n)}a_2^{p_2(n)}$
 are polynomial sequences in $H$.

Let $\pi\colon H/\Delta\to H/([H_0,H_0]\Delta)$ be the natural
projection.   By Theorem~\ref{T:L} it suffices to  show that for
every $x\in X$ we have
$\overline{\{h_2(m,n)\pi(\tilde{x})\}}_{m,n\in\N}=
H/([H_0,H_0]\Delta)$. Since the sequence
$\{h_1(m,n,r)\pi(\tilde{x})\}_{m,n,r\in\N}$ is dense in
$H/([H_0,H_0]\Delta)$ it is enough to show that for every $x\in X$
we have
\begin{equation}\label{mem}
\overline{\{h_2(m,n)\pi(\tilde{x})\}}_{m,n\in\N}=
\overline{\{h_1(m,n,r)\pi(\tilde{x})\}}_{m,n,r\in\N}.
\end{equation}

We first obtain some information about the quotient
$H/([H_0,H_0]\Delta)$. We  claim that $[H_0,H_0]=[G_0,G_0]^4$.
 The
$\subset$ inclusion is obvious. To establish the other inclusion
first notice that for $g\in G_0$ elements of the form $(g,g,g,g)$,
$(e,g^{k_1},e,g^{l_1})$, and $(e,e,g^{k_2},g^{l_2})$ belong to
$H_0$. Taking commutators of these elements and using the fact
that the group $G_0$ is divisible, we easily get that
$[G_0,G_0]^4\subset [H_0,H_0]$, proving the  claim.
Hence, the quotient space $ H/([H_0,H_0]\Delta)$ can be identified
with a subset of $Z^4$
 where $Z=G/([G_0,G_0]\Gamma)$, and so we can consider both sets in \eqref{mem} as subsets of
$Z^4$.

Since $Z$ is connected and the connected component of the group
$G/[G_0,G_0]$ is abelian,
 Theorem~\ref{T:FK} applies. So we can assume that $Z=\T^d$
for some $d\in \N$, $\pi(\tilde{x})$ is represented by some
$(x,x,x,x)\in \T^{4d}$, and the nilrotations by $(a,e,e,e)$,
$(e,a,e,e)$, $(e,e,a,e)$, and $(e,e,e,a)$, are represented by the
transformations $S\times\text{id}\times\text{id}\times \text{id}$,
$\text{id}\times S\times\text{id}\times \text{id}$,
$\text{id}\times\text{id}\times S\times \text{id}$,
$\text{id}\times\text{id}\times \text{id}\times S$,
where $S$ is an ergodic nilpotent affine transformation of $\T^d$.
We have thus reduced our problem to showing that for every ergodic
nilpotent affine transformation $S$ acting on $X=\T^d$, linearly
independent integer polynomials $p_1,p_2$, and every $x\in \T^d$,
the sequences
$$
\{(S^mx,S^{m+k_1n}x,S^{m+k_2r}x,S^{m+l_1n+l_2r}x)\}_{m,n,r\in \N}
$$
and
$$
\{(S^mx,S^{m+k_1p_1(n)}x,S^{m+k_2p_2(n)}x,S^{m+l_1p_1(n)+l_2p_2(n)}x)\}_{m,n\in
\N}
$$
have the same closure. Since $S$ is uniquely ergodic the sequence
$\{S^mx\}_{m\in\N}$ is dense in $X$ for every $x\in X$.  So it
suffices to show that the sets $\mathcal{O}(P,\Delta_{X^3},S)$ and
$\mathcal{O}(Q,\Delta_{X^3},S)$ have the same closure,  where $Q$
is the family of 2-variable polynomials $\{k_1n,k_2r,l_1n+l_2r\}$.
This in turn will follow if we show that the averages \eqref{b1}
and \eqref{b2} have the same limit as $N-M\to\infty$ in the
special case where the transformation $T$ is equal to $S$. Since
$W(P)=2$, by Proposition~\ref{P:weyl3} the characteristic factor
for the averages \eqref{b1} when $T=S$ is the Kronecker factor. By
Lemma~\ref{L:w2'} the Kronecker factor  is also characteristic for
the averages \eqref{b2}, so it suffices to check the identity for
group rotations. This can be easily verified for characters and
then for general bounded functions by approximating them in $L^2$
by finite linear combinations of characters, thus completing the
proof.
\end{proof}

\subsection{Characteristic factors and limit formulas}
We are now ready to  prove Theorem~B. The argument is rather
lengthy  so we refer the reader to the Introduction for a brief
sketch. Notice that by Proposition~\ref{T:weyl} the cases $(i)$,
$(ii)$, $(iii)$ of Theorem~B correspond to the cases where the
polynomial family has  Weyl complexity $1$, $3$, $2$
correspondingly. We deal with each one separately.

\subsubsection{{\bf Weyl complexity $1$}}
\emph{Characteristic factor}: We can assume that $p_i(0)=0$ for
$i=1,2,3$. If the polynomials are linearly independent it was
shown in \cite{FK2} that the rational Kronecker factor
$\mathcal{K}_{rat}$ is characteristic  for $L^2$-convergence of
the averages $(P)$.

\emph{Limit formula}: In the case where the system is totally
ergodic the factor $\mathcal{K}_{rat}$ is trivial, hence for every
$f_1,f_2,f_3\in L^\infty(\mu)$ we have
\begin{equation}\label{E:form1}
\lim_{N-M\to\infty} \frac{1}{N-M}\sum_{n=M}^{N-1}
f_1(T^{p_1(n)}x)\cdot f_2(T^{p_2(n)}x) \cdot f_3(T^{p_3(n)}x)=
\int f_1\ d\mu\cdot \int f_2 \ d\mu \cdot \int f_3 \ d\mu,
\end{equation}
where the limit is taken in $L^2(\mu)$.


\subsubsection{{\bf Weyl complexity $2$}}
\emph{Characteristic factor}: The collection
   of $3$-term polynomial families of Weyl complexity $2$
   is easily shown to be eligible, so
   by Proposition~\ref{P:te} we can
   assume that the system is totally ergodic.
It follows from  Proposition~\ref{T:weyl} that the polynomials $\{\tilde{p}_1,\tilde{p}_2,\tilde{p}_3\}$
are linearly dependent. Hence,
 $$
 (p_1,p_2,p_3)=(k_1q_1+c_1, k_2q_2+c_2, l_1q_1+l_2q_2+c_3)
 $$
  for some linearly independent integer polynomials  $q_1,q_2$ with zero constant term
  and
 $k_1,k_2,l_1,l_2,c_1,c_2,c_3\in\Z$.
Combining Lemmas~\ref{L:w2} and \ref{L:w2'}, we get that for
totally ergodic systems the Kronecker factor $\mathcal{K}$ is
characteristic for $L^2$-convergence of the averages $(P)$.

It can be easily seen that for polynomial families of Weyl
complexity $2$ every  characteristic factor (thought of as a
subalgebra of functions) for the averages $(P)$ contains all the
eigenfunctions of the system, and as a result it  contains the
Kronecker factor. Hence, for ergodic systems the Kronecker factor
is the smallest characteristic factor.

\emph{Limit formula}:  We now compute the limit of the
corresponding ergodic  averages  $(P)$ for totally ergodic
systems. We can assume that $c_i=0$ for $i=1,2,3$.
 After replacing all three functions with their projection to the
Kronecker factor $\mathcal{K}$ we can assume that
$\mathcal{X}=\mathcal{K}$. Every Kronecker system is an inverse
limit of $1$-step nilsystems so we can assume that our system  is
a totally ergodic rotation on a compact abelian Lie group $G$ with
the Haar measure $m$. Moreover, by Proposition~\ref{P:te} the
group $G$ has to be connected, so $G=\T^d$ for some nonnegative
integer $d$. In this case it is easy to check that for every
$f_1,f_2,f_3\in L^\infty(\mu)$ we have that
\begin{multline}\label{E:gfg}
\lim_{N-M\to\infty} \frac{1}{N-M}\sum_{n=M}^{N-1}
f_1(T^{k_1q_1(n)}t)\cdot
f_2(T^{k_2q_2(n)}t) \cdot f_3(T^{l_1q_1(n)+l_2q_2(n)}t)=\\
\int_{\T^{2d}} f_1(t+k_1t_1)\cdot f_2(t+k_2t_2)\cdot
f_3(t+l_1t_1+l_2t_2)\ dm(t_1) \ dm(t_2)
\end{multline}
for a.e. $t\in \T^d$.

\subsubsection{{\bf Weyl complexity $3$}} \label{SSS:weyl3}
 By Proposition~\ref{T:weyl} the polynomial triple $(p_1,p_2,p_3)$ either
has the form
$(a) \ (lp+c_1,mp+c_2,kp^2+rp+c_3)$,   or
$(b)\ (kp^2+lp+c_1,kp^2+mp+c_2,kp^2+rp+c_3)$, for some integer
polynomial $p$, and  $k,l,m,r,c_1,c_2,c_3\in\Z$. We consider the
following three cases:

 {\bf Case 1}: The family of essentially distinct polynomials has the form
 $(a)$ with $k=0$. This case is covered by Theorem~A. The smallest
 characteristic factor is $\mathcal{Z}_2$.
 To find a limit
 formula in the totally ergodic case, first using standard deductions
 we can assume that the
 system is a totally ergodic $2$-step nilsystem. In this case
 Theorem~\ref{T:Z}  gives a formula for the limit.

{\bf Case 2}: The family of essentially distinct polynomials has
the form $(a)$ with $k\neq 0$. We first deal with the case
$p(n)=n$ and then reduce the case of a general polynomial $p(n)$
to this one.

  \emph{Characteristic factor for $p(n)=n$}:
  Since the collection
   of polynomial families of the form $(a)$ with $k\neq 0$ is eligible,
   by Proposition~\ref{P:te} we can
   assume that the system is totally ergodic. Furthermore, we can
   assume that $c_i=0$ for $i=1,2,3$.
    We first claim that if $f_3\in \mathcal{K}^{\bot}$ then the
averages
\begin{equation}\label{E:av1}
 \frac{1}{N-M}\sum_{n=M}^{N-1} T^{ln}f_1\cdot T^{mn}f_2\cdot
T^{kn^2+rn}f_3
\end{equation}
 converge to zero in $L^2(\mu)$ as $N-M\to\infty$.
 We apply the Hilbert space  Van
der Corput Lemma (\cite{Be5})\footnote{Let $\{x_n\}_{n\in \N}$ be
a bounded sequence in a Hilbert space.
 If for every $m\in \N$ one has
 $
  \lim_{N-M\to\infty}\frac{1}{N-M} \sum_{n=M}^{N-1} \langle
x_{n+m},x_n\rangle=0,
 $
then $ \lim_{N-M\to\infty}\norm{\frac{1}{N-M} \sum_{n=M}^{N-1}
x_n}=0. $}
  for the sequence of functions
$$
a_n(x)=f_1(T^{ln}x)\cdot f_2(T^{mn}x)\cdot f_3(T^{kn^2+rn}x).
$$
It suffices to show that for every $h\in \N$ we have
$$
\lim_{N-M\to\infty}\frac{1}{N-M}\sum_{n=M}^{N-1}\langle
a_{n+h},a_n\rangle=0,
$$
or equivalently that the average
\begin{equation}\label{opa}
 \frac{1}{N-M}\sum_{n=M}^{N-1} T^{(m-l)n}(T^{mh}f_2\cdot \bar{f}_2)\cdot
 T^{kn^2+(2kh+r)n}(T^{kh^2+h}f_3)\cdot T^{kn^2+rn}\bar{f}_3
\end{equation}
converges to zero in $L^2(\mu)$ as $N-M\to\infty$. Using
Proposition~\ref{T:weyl} it is easy to check  that for all
$h,k,l,m,r\in\Z$ with  $h,k,l,m\neq 0$ and $l\neq m$ we have
$$
W((m-l)n, kn^2+(2kh+r)n,kn^2+rn)= 2.
$$
Hence, as shown in the Weyl complexity $2$ case, the
characteristic factor for the  ergodic averages~\eqref{opa} is the
Kronecker factor. This proves the claim.

Next we claim that if $f_1$ or $f_2 \in\mathcal{A}_2^\bot$  then
the averages \eqref{E:av1} converge to zero in $L^2(\mu)$ as
$N-M\to\infty$. We prove this for $f_1$, the argument is similar
for $f_2$. As we have shown, we can replace $f_3$ with
$\E(f_3|\mathcal{K})$ without changing the limit of the averages
\eqref{E:av1}. Moreover, after approximating $\E(f_3|\mathcal{K})$
by a linear combination of eigenfunctions, and using linearity, we
can assume that $f_3$ is either constant, or  a
$\lambda$-eigenfunction where $\lambda=e(\alpha)$ for some
$\alpha\in (0,1)$. Moreover, since the system is totally ergodic
$\alpha$ is irrational.  If $f_3$ is constant the claim follows
from a classical result of Furstenberg~\cite{Fu1}. If not, the
average \eqref{E:av1} is equal to $f_3$ times the average
\begin{equation}\label{E:reduced1}
\frac{1}{N-M}\sum_{n=M}^{N-1} T^{ln}f_1\cdot T^{mn}f_2\cdot
e((kn^2+rn)\alpha).
\end{equation}
 A simple computation shows that there exist characters
$\chi_1,\chi_2\colon \mathbb{T}^2\to\mathbb{C}$ such that
$$
\chi_1(R^{ln}(t_1,t_2))\cdot \chi_2(R^{mn}(t_1,t_2))=
e((kn^2+rn)\alpha)
$$
holds for every $n\in\N$, where
$R\colon\mathbb{T}^2\to\mathbb{T}^2$ is defined by
$$
R(t_1,t_2)=(t_1+\beta,t_2+2t_1+\beta)
$$
and $\beta$ is some appropriately chosen rational multiple of
$\alpha$. Consider the product system
$(X\times\mathbb{T}^2,\mu\times m, S=T\times R)$, where $m$ is the
Haar measure on $\T^2$, and let $h_1=f_1\cdot \chi_1$,
$h_2=f_2\cdot \chi_2$. Then the average \eqref{E:reduced1} takes
the form
\begin{equation}\label{E:reduced2}
\frac{1}{N-M}\sum_{n=M}^{N-1} S^{ln}h_1\cdot S^{mn}h_2.
\end{equation}
Let $S_t, t\in [0,1],$ be the ergodic components of  $S$. We will
show that if  $f_1\in \mathcal{A}_2(T)^\bot$ then  $h_1(x)\in
\mathcal{K}(S_t)^\bot$ for a.e. $t$. As it is well known, this
would follow if we show that for a.e. $\tilde{x}\in X\times \T^2$
we have
$$
\lim_{N\to\infty} \sup_{s\in [0,1)} \Big|\frac{1}{N}
\sum_{n=0}^{N-1} h_1(S^n\tilde{x})\cdot e(ns)\Big|=0,
$$
or equivalently that
$$
\lim_{N\to\infty} \sup_{s\in [0,1)} \Big|\frac{1}{N}
\sum_{n=0}^{N-1} f_1(T^nx)\cdot e(ns+n^2\gamma)\Big|=0
$$
for a.e. $x\in X$, where $\gamma$ is some integer multiple of
$\beta$. Since $f_1\in \mathcal{A}_2(T)^\bot$ and $T$ is totally
ergodic this follows from \cite{Fr}.  Hence,  $h_1(x)\in
\mathcal{K}(S_t)^\bot$ for a.e. $t$. From \cite{Fu1} we know that
for distinct nonzero integers $l,m$ the Kronecker factor is
characteristic for the ergodic averages associated to the family
$\{ln,mn\}$. So an ergodic decomposition argument gives that the
average in \eqref{E:reduced2} converges to zero in $L^2(\mu)$ as
$N-M\to\infty$, proving the claim. This shows that the factor
$\mathcal{A}_2$ is characteristic for the averages
 \eqref{E:av1}.

\emph{Limit formula for $p(n)=n$}: We now compute the limit of the
averages \eqref{E:av1} for totally ergodic systems. We can assume
that $c_i=0$ for $i=1,2,3$. Since $\mathcal{A}_2$ is a factor of
$\mathcal{Z}_2$, and $\mathcal{Z}_2$ is an inverse limit of
$2$-step nilsystems, using an approximation argument we can assume
that our system is a totally ergodic $2$-step nilsystem that
coincides with its $2$-step affine factor $\mathcal{A}_2$. In this
case the system is isomorphic to a $2$-step nilpotent affine
transformation on a connected compact abelian group $G$ (see
\cite{A}). Furthermore,  our system is a nilsystem, so the group
$G$ has to be Lie. Hence, we can assume that $G$ is a finite
dimensional torus. In this case, the  evaluation of the limit is a
straightforward computation, which is done (for general $k$-step
affine systems) in \cite{BLL} (or \cite{L4}).
Instead of reproducing this
rather complicated formula let
 us illustrate how the limit is computed in a simple case.
Suppose that $T\colon \T^2\to \T^2$ is given by
$$
T(t_1,t_2)=(t_1+\alpha,t_2+b
t_1)
$$
where $\alpha$ is irrational and $b$ is a nonzero integer. We find
by direct computation that for a.e. $(t_1,t_2)\in\T^2$ we have
\begin{multline}\label{E:dfd}
\lim_{N-M\to\infty} \frac{1}{N-M}\sum_{n=M}^{N-1}
f_1(T^{ln}(t_1,t_2))\cdot
f_2(T^{mn}(t_1,t_2)) \cdot f_3(T^{kn^2+rn}(t_1,t_2))=\\
\int_{\T^{4}} f_1(t_1+lx_1, t_2+ly_1+l^2y_2)\cdot
f_2(t_1+mx_1,t_2+my_1+m^2y_2)\cdot f_3(t_1+ky_2+rx_1,y_3)\ dx_1 \
dy_1 \ dy_2\ dy_3.
\end{multline}
It immediately
follows from this formula (and more generally from the formulas in \cite{BLL} or \cite{L4})
that for almost every $x\in G$
the set $H_x$ of \eqref{E:H_x} is connected.

\emph{Connectedness for $p(n)=n$}: In order to deal with the case
of a general polynomial $p(n)$ we will apply
Proposition~\ref{L:p(n)} which allows as to make the substitution
$n\to p(n)$ when computing the orbit closure of a polynomial
sequence with connected closure. We now verify that the
connectedness assumption is satisfied, i.e. that
 for every  totally ergodic
nilsystem $(X=G/\Gamma,\G/\Gamma,m,T_a)$ the set
\begin{equation}\label{E:H_x}
H_x=\overline{\{(a^{ln}x,a^{mn}x,a^{kn^2+rn}x)\}}_{n\in\N}
\end{equation}
 is connected for a.e. $x\in X$.
 By
Theorem~\ref{T:L} we have
\begin{equation}\label{E:hjh}
\lim_{N-M\to\infty} \frac{1}{N-M}\sum_{n=M}^{N-1}
f_1(a^{ln}x)\cdot f_2(a^{mn}x)\cdot f_3(a^{kn^2+rn}x)= \int_{H_x}
f_1\otimes f_2\otimes f_3 \ dm_{H_x}
\end{equation}
for a.e. $x\in X$. Since the factor $\mathcal{A}_2$ is
characteristic for convergence of the averages in \eqref{E:hjh},
we can replace every function by its projection to $\mathcal{A}_2$
which by Proposition~\ref{P:nilsystem} is
$Z=G/([G_0,G_0]G_3\Gamma)$. This shows that the set $H_x$ factors
through $Z^3$.   Furthermore, we know that $T_a$ acting on $Z$ is
topologically conjugate to a $2$-step nilpotent  affine
transformation on some finite dimensional torus $\T^d$. As we
mentioned before, we can compute explicitly the limit in this case
and derive that for a.e. $x\in X$ the projection of $\pi(H_x)$ of
$H_x$ onto $Z^3$ is connected. It follows that  the set $H_x$ is a product
of the connected set $\pi(H_x)$ and the connected nilmanifold
$([G_0,G_0]G_3)/\Gamma'$, where $\Gamma'= \Gamma \cap
([G_0,G_0]G_3)$. Hence, for a.e. $x\in X$ the set $H_{x}$ is
connected.

\emph{General case}: To deal with the general case notice that all
the previous results carry through once we show that for totally
ergodic systems the $L^2$-limit of the averages in
\eqref{E:reduced1} remains the same if we replace $n$ with any
nonconstant polynomial $p(n)$. Using  Theorem~\ref{T:HKL} and an
approximation argument, it suffices to verify  that this is the
case for totally ergodic nilsystems. By Theorem~\ref{T:L} we can
further reduce this to showing that if
$(X=G/\Gamma,\mathcal{G}/\Gamma,m,T_a)$ is a nilsystem with $X$
connected, then for almost every $x\in X$ the sequences $
\{(a^{ln}x, a^{mn}x, a^{kn^2+rn}x)\}_{n\in \N} $ and $
\{(a^{lp(n)}x, a^{mp(n)}x, a^{kp(n)^2+rp(n)}x)\}_{n\in \N} $ have
the same closure. We previously showed that for a.e. $x\in X$ the
closure of the first sequence is connected. Hence,
Proposition~\ref{L:p(n)} applies and proves the claim.

It can be easily seen that for polynomial families of the form
$(a)$ with $k\neq 0$  every characteristic factor (thought of as a
subalgebra of functions) for the averages $(P)$ contains all the
functions in $\mathcal{E}_2$ (defined in
Section~\ref{SS:ergodic}), and as a result it contains the factor
$\mathcal{A}_2$. Hence, for ergodic systems the factor
$\mathcal{A}_2$ is the smallest characteristic factor.

{\bf Case 3}:  The family of essentially distinct polynomials has
the form
 $(b)$ with $k\neq 0$.

\emph{Characteristic factor}:
     It suffices to show that if
$f_i \in\mathcal{A}_2^\bot$ for some $i=1,2,3$ then the averages
\begin{equation}\label{E:av3}
 \frac{1}{N-M}\sum_{n=M}^{N-1} T^{kp(n)^2+lp(n)}f_1\cdot T^{kp(n)^2+mp(n)}f_2\cdot
T^{kp(n)^2+rp(n)}f_3
\end{equation} converge to zero in
$L^2(\mu)$ as $N-M\to\infty$. We show this for $i=1$, the argument
is similar for $i=2,3$. This time applying Van der Corput's lemma
doesn't help. Instead, we notice that since the limit in
$L^2(\mu)$ as $N-M\to\infty$ of the averages  \eqref{E:av3}
exists, it suffices to show that if $f_1 \in\mathcal{A}_2^\bot$
then for every $f_0\in L^\infty(\mu)$ we have
\begin{equation*}
\lim_{N-M\to\infty} \frac{1}{N-M}\sum_{n=M}^{N-1} \int f_0\cdot
T^{kp(n)^2+lp(n)}f_1\cdot T^{kp(n)^2+mp(n)}f_2\cdot
T^{kp(n)^2+rp(n)}f_3\ d\mu=0.
\end{equation*}
Equivalently we need to show that
\begin{equation*}
\lim_{N-M\to\infty} \frac{1}{N-M}\sum_{n=M}^{N-1} \int f_3\cdot
T^{(l-r)p(n)}f_1\cdot T^{(m-r)p(n)}f_2\cdot T^{-kp(n)^2-rp(n)}f_0\
d\mu=0
\end{equation*}
for every $f_0\in L^\infty(\mu)$, which is true by Case $2$.

An argument analogous to the one explained in Case $2$ shows that
for ergodic systems the factor $\mathcal{A}_2$  is the smallest
characteristic factor for the averages \eqref{E:av3}. Also a limit
formula goes along the lines of Case $2$.

\section{Applications in combinatorics}\label{S:applications}
In this section we are going to derive several  combinatorial
implications of our results in ergodic theory. Our starting point
will always be the Correspondence Principle of Furstenberg that
enables us to translate statements in combinatorics to statements
in ergodic theory. We mention a slight modification of this
principle due to Lesigne (see \cite{BHK}) that allows us to work
with ergodic systems (this is crucial for Theorem~C'):
\begin{Correspondence1}
For every $\Lambda\subset \N$  there exists an invertible ergodic
system $(X,\X,\mu,T)$ and $A\in\mathcal{X}$ with
$\mu(A)=d^*(\Lambda)$
 and such that
$$
d^*(\Lambda \cap (\Lambda+n_1)\cap\ldots\cap (\Lambda+n_k))\geq
\mu(A\cap T^{n_1}A\cap\cdots \cap T^{n_k}A),
$$
for all $k\in\mathbb{N}$ and integers $n_1,\ldots,n_k$.
\end{Correspondence1}

\subsection{Sets of multiple recurrence}
We will prove Theorem~D.
\begin{proof}[Proof of Theorem~D]
Suppose that $p(n)$ is an integer polynomial that satisfies the
assumptions of the theorem. Using Furstenberg's Correspondence
Principle it suffices to show that if $f\in L^\infty(\mu)$ is
nonnegative and not a.e. zero then
\begin{equation}\label{positive}
\lim_{N\to\infty} \frac{1}{N} \sum_{n=0}^{N-1} \int f\cdot
T^{p(n)}f\cdot \ldots \cdot T^{kp(n)}f\ d\mu>0.
\end{equation}
Using  an ergodic decomposition argument we can assume that the
system is ergodic and by  Theorem~\ref{T:HKL} we can reduce the
problem to showing \eqref{positive} in the case where the system
is an inverse limit of nilsystems. Moreover, an argument
completely analogous to that of Lemma $3.2$ in \cite{FuK} shows
that the positiveness property \eqref{positive} is preserved by
inverse limits.
Hence, we can further assume that the system is an ergodic
nilsystem. In this case by Proposition~\ref{P:te} there exists an
$r\in \N$ such that the  ergodic components of $T^r$ are totally
ergodic. By our assumption there exists an $n_0\in \N$ such that
$p(n_0)\equiv 0 \pmod{r}$. Then $p(rn+n_0)=rq(n)$ for some integer
polynomial $q$ and the limit in \eqref{positive} is greater or
equal than $1/r$ times
$$
\lim_{N\to\infty} \frac{1}{N} \sum_{n=0}^{N-1} \int f\cdot
T^{rq(n)}f\cdot \ldots \cdot T^{krq(n)}f\ d\mu.
$$
Using Theorem~A for the ergodic components of $T^r$ we get that
this last limit  equals
$$
\lim_{N\to\infty} \frac{1}{N} \sum_{n=0}^{N-1} \int f\cdot
T^{n}f\cdot \ldots \cdot T^{kn}f\ d\mu $$ which  is positive by
\cite{Fu1}.
\end{proof}

\subsection{A bad set for recurrence with good powers}
We will prove Theorem~E. It will be a consequence of the
Polynomial  Szemer\'edi Theorem and the following multiple ergodic
theorem:

\begin{proposition}\label{P:squares}
 Let $(X,\mathcal{X},\mu,T)$ be an invertible
   system, $h\colon
\mathbb{T}\to\mathbb{C}$ be Riemann integrable,
 $f_0,\ldots,f_k\in L^\infty(\mu)$,  and $\beta$ be an
irrational number. Then for every integer polynomial $p$ with
$\deg{p}>1$ we have
\begin{gather}\label{E:main2}
\lim_{N\to\infty} \frac{1}{N} \sum_{n=0}^{N-1}  h(n\beta)\cdot
\int f_0\cdot T^{p(n)}f_1\cdot \ldots \cdot T^{kp(n)}f_k\ d\mu=\\
\notag \int h\ dt\cdot  \lim_{N\to\infty} \frac{1}{N}
\sum_{n=0}^{N-1} \int f_0\cdot T^{p(n)}f_1\cdot \ldots \cdot
T^{kp(n)}f_k \ d\mu.
\end{gather}
\end{proposition}
\begin{proof}
  Using  an ergodic decomposition argument it suffices to check
  \eqref{E:main2}  when the system is ergodic.
  In \cite{HK} it is shown that for every $k\in\N$ if $\mathbb{E}(f
  |\mathcal{Z}_m(T))=0$ then
  $\mathbb{E}(f\otimes f|\mathcal{Z}_{m-1}((T\times T)_t))$  for
  a.e.  $t$, where $(T\times T)_t$, $t\in [0,1]$, denote
  the ergodic components of $T\times T$. Keeping this in mind, and
  applying    Theorem~\ref{T:HKL}
  for the ergodic components of the
  product system $(X\times X,\mathcal{X}\times\mathcal{X}, \mu\times
  \mu,T\times T)$, we get that there exists an $m\in\N$ such that
  if $E(f_i|\mathcal{Z}_m)=0$ for some $i=0,\ldots,k$ then
  $$
 \text{D-}\!\!\lim_{n\to\infty}\int f_0\cdot
T^{p_1(n)}f_1\cdot \ldots \cdot T^{p_k(n)}f_k \ d\mu=0.
$$
So \eqref{E:main2} is obvious when $E(f_i|\mathcal{Z}_m)=0$ for
some $i=0,\ldots,k$, since then both limits   are zero. We can
therefore assume that $f_i\in \mathcal{Z}_m$ for $i=0,\ldots,k$.
Since the factor $\mathcal{Z}_m(T)$ is an inverse limit of
nilsystems, a standard approximation argument shows that it
suffices to check \eqref{E:main2} when the system is an ergodic
nilsystem, say $(X=G/\Gamma,\mathcal{G}/\Gamma,m,T_a)$.
In this case, equation \eqref{E:main2}  follows if we show that
for $f_1,\ldots,f_k\in L^\infty(\mu)$ we have for a.e. $x\in X$
that
\begin{gather}\label{E:squares}
\lim_{N\to\infty} \frac{1}{N} \sum_{n=0}^{N-1}  h(n\beta)\cdot
f_1(a^{p(n)}x)\cdot \ldots \cdot f_k(a^{kp(n)}x)= \\ \notag \int
h(t)\ dt\cdot \ \lim_{N\to\infty} \frac{1}{N} \sum_{n=0}^{N-1}
f_1(a^{p(n)}x)\cdot\ldots\cdot f_k(a^{kp(n)}x).
\end{gather}
By Theorem~\ref{T:L} it suffices to show that for a.e. $x\in X$ we
have
\begin{equation}\label{E:uiu}
\overline{\{(n\beta, a^{p(n)}x,\ldots, a^{kp(n)}x)\}}_{n\in\N}=
\T\times \overline{\{(a^{p(n)}x,\ldots, a^{kp(n)}x)\}}_{n\in\N}.
\end{equation}
Since $\deg{p}>1$ this follows from Lemma~\ref{L:squares}.
\end{proof}

\begin{proof}[Proof of the Theorem~E]
We will show that the set $S=\big\{n\in\N\colon \{n\sqrt{2}\}\in
[1/4,3/4]\big\}$ has the advertised property. Clearly  $S$ is not
good for single recurrence since it is not good for
 recurrence for the rotation by $\sqrt{2}$ on $\mathbb{T}$.

We will show that $p(S)$ is a set of multiple recurrence whenever
$p$ is an integer polynomial with $\deg{p}>1$. So let
$(X,\mathcal{X},\mu,T)$ be an invertible system and
$A\in\mathcal{X}$ with $\mu(A)>0$. We apply
Proposition~\ref{P:squares} for $f_i={\bf 1}_A$, $i=0,1,\ldots,k$,
$h={\bf 1}_{[1/4,3/4]}$, and $\beta=\sqrt{2}$. We get
\begin{gather}
\lim_{N\to\infty} \frac{1}{N} \sum_{n\in S\cap [1,N]} \mu(A\cap
T^{-p(n)}A\cap\cdots\cap T^{-kp(n)}A)= \\ \frac{1}{2}\cdot \notag
\lim_{N\to\infty} \frac{1}{N} \sum_{n=0}^{N-1} \mu(A\cap
T^{-p(n)}A\cap\cdots\cap T^{-kp(n)}A).
\end{gather}
The last limit is positive by Theorem~\ref{T:PSzemeredi}, showing
that $p(S)$ is a set of multiple recurrence.
\end{proof}
\subsection{Universal families of three polynomials}
We prove Theorem~F.
\begin{proof}[Proof of Theorem~F]
We can assume that  the polynomials $p_1, p_2, p_3$ are
essentially distinct.
 We claim that under the assumptions of the
theorem, if $f\in L^\infty(\mu)$ is nonnegative and not a.e. zero
then
\begin{equation}\label{positive2}
\lim_{N\to\infty} \frac{1}{N} \sum_{n=0}^{N-1} \int f\cdot
T^{p_1(n)}f\cdot T^{p_2(n)}f\cdot T^{p_3(n)}f d\mu>0.
\end{equation}
An argument analogous to the one used in the beginning of the
proof of Theorem~D allows us to reduce the problem to showing
\eqref{positive2} in the case where the system is an ergodic
nilsystem, say  $(X=G/\Gamma,\mathcal{G}/\Gamma,m,T_a)$.

{\bf Weyl complexity $1$.} If the polynomials
$p_1-p(0),p_2-p_2(0),p_3-p_3(0)$ are linearly independent we have
from Theorem~B that the factor $\cK_{rat}$ is characteristic for
the averages in \eqref{positive2}, hence we can assume that
$\mathcal{X}=\cK_{rat}$. Since our system is a nilsystem we have
$\cK_{rat}=\cK_r$ for some $r\in\N$. By our assumption there
exists $n_0\in\N$ such that $p_i(n_0)\equiv 0 \pmod{r}$ for
$i=1,2,3$. Then $p_i(rn+n_0)=rp_i'(n)$ for some  integer
polynomials $p_i'$ for $i=1,2,3$. Hence,  whenever $n\equiv n_0
\pmod{r}$ we have $T^{p_i(n)}=\text{id}$ for $i=1,2,3$, and so the
integral in \eqref{positive2} is equal to $\int f^4\ d\mu>0$. The
result follows.

 {\bf Weyl complexity $2$.} We start with some reductions on the
polynomial family. We have that $(p_1,p_2,p_3)= (k_1q_1+c_1,
k_2q_2+c_2, l_1q_1+l_2q_2+c_3)$ for some  linearly independent
integer polynomials $q_1,q_2$ and
$k_1,k_2,l_1,l_2,c_1,c_2,c_3\in\Z$. Since $p_i(n)\equiv 0
\pmod{k_i}$ has a solution for $i=1,2$, we get that $c_1=k_1c_1'$,
$c_2=k_2c_2'$. So we are reduced to the case where the polynomial
family has the form $(k_1q_1, k_2q_2, l_1q_1+l_2q_2+c_3)$. If
$c_3\neq 0$ we can choose an $r\in\N$ that is relatively prime to
the integers $k_1,k_2,c_3$. Then the system of equations
$p_i(n)\equiv 0 \pmod{r}$, $i=1,2,3$, does not have a solution,
contrary to our assumption. Hence, $c_3=0$. So we can assume that
$$
(p_1,p_2,p_3)=(k_1q_1, k_2q_2, l_1q_1+l_2q_2)
$$
 for some linearly
independent integer polynomials $q_1,q_2$ and $k_1,k_2,l_1,l_2\in
\Z$.

 By Proposition~\ref{basic} there exists an $r\in \N$ such that
the ergodic components of $T_a^r$ are totally ergodic. By our
assumption there exists $n_0\in\N$ such that $q_i(n_0)\equiv 0
\pmod{r}$ for $i=1,2$. Then $q_i(rn+n_0)=rq_i'(n)$ for some
linearly independent integer polynomials $q_1', q_2'$, and
 the average in \eqref{positive2} is greater or equal than $1/r$
times
\begin{equation}\label{positive3}
\liminf_{N\to\infty} \frac{1}{N} \sum_{n=0}^{N-1} \int f\cdot
T_a^{rk_1q_1'(n)}f\cdot T_a^{rk_2q_2'(n)}f\cdot
T_a^{r(l_1q_1'(n)+l_2q_2'(n))}f \ d\mu.
\end{equation}
Working with the (totally ergodic)  ergodic components of $T_a^r$
and the polynomial family $\{k_1q_1',k_2q_2',l_1q_1'+l_2q_2'\}$
(which also has Weyl complexity $2$) we get from Lemma~\ref{L:w2}
that the limit in \eqref{positive3} is equal to
$$
\lim_{N\to\infty} \frac{1}{N^2} \sum_{n,r=1}^N\int f \cdot
T_a^{k_1n}f\cdot T_a^{k_2r}f\cdot T_a^{l_1n+l_2r}f\ d\mu
$$
which is easily shown to be positive.

  {\bf Weyl complexity $3$.}
 The argument is similar to the one used in the previous case so we just sketch
 the main steps.  By Proposition~\ref{T:weyl} some permutation of the polynomials
  $p_1,p_2,p_3$ either
have the form $(a) \ (lp+c_1,mp+c_2,kp^2+rp+c_3)$, or the form
$(b)\ (kp^2+lp+c_1,kp^2+mp+c_2,kp^2+rp+c_3)$, for some integer
polynomial $p$, and some $k,l,m,r,c_1,c_2,c_3\in\Z$. Arguing as in
the Weyl complexity $2$ case, we can assume that $c_i=0$ for
$i=1,2,3$ and the system is totally ergodic. We consider the
following three cases:

In the case  $(a)$ with $k=0$ we get from Theorem~A that the limit
in \eqref{positive2} is equal to
$$
\lim_{N\to\infty} \frac{1}{N}\sum_{n=0}^{N-1} \int f\cdot
T^{k_1n}f \cdot T^{k_2n}f\cdot T^{k_3n}f \ d\mu,
$$
which is positive by \cite{Fu1}.

In the case  $(a)$ with $k\neq 0$  we showed in the proof of part
$(ii)$ of Theorem~B that the limit in \eqref{positive2} is equal
to
$$
\lim_{N\to\infty} \frac{1}{N}\sum_{n=0}^{N-1} \int f\cdot T^{ln}f
\cdot T^{mn}f\cdot T^{kn^2+rn}f \ d\mu,
$$
which is positive by Theorem~\ref{T:PSzemeredi}.

To deal with  the case $(b)$ with $k\neq 0$  we use the identity
$$
\int f \cdot T^{p_1(n)}f\cdot T^{p_2(n)}f\cdot T^{p_3(n)}f \ d\mu=
\int T^{-p_3(n)}f \cdot T^{p_1(n)-p_3(n)}f \cdot
T^{p_2(n)-p_3(n)}f \cdot f \ d\mu
$$
which  allows us to  show to reduce case $(b)$ with $k\neq 0$ to
case $(a)$ with $k\neq 0$ that we previously handled. This
completes the proof.
\end{proof}

\subsection{Positive results for lower bounds}
 The proof of Theorem~C' is an immediate consequence
 of  Theorem~C
 and  Furstenberg's Correspondence Principle.
  So it remains to prove Theorem~C.
\begin{proof}[Proof of Theorem~C]
The proof for the case of two polynomials goes along the lines of
the case of three polynomials with Weyl complexity $\leq 2$ and so
we omit it.

So let $\{p_1,p_2,p_3\}$ be a family of essentially distinct
integer polynomials that is not equal to any of  the exceptional
forms mentioned in  Theorem~C. Then by Proposition~\ref{T:weyl}
the polynomial family either has Weyl complexity $\leq 2$, or some
permutation of the polynomials has the form $\{kp,lp,(k+l)p\}$,
for some integer polynomial $p$ with $p(0)=0$ and $k,l\in\Z$. So
we have to deal with the following two cases:

{\bf Case 1.} Suppose that $W(p_1,p_2,p_3)\leq 2$. If
$W(p_1,p_2,p_3)= 1$ the polynomials are linearly independent and
the result follows from \cite{FK2}. If $W(p_1,p_2,p_3)=2 $ we can
assume that
 $(p_1,p_2,p_3)= (k_1q_1, k_2q_2, l_1q_1+l_2q_2)$, where $q_1,q_2$ are some
linearly independent integer polynomials and  $k_1,k_2,l_1,l_2\in
\Z$.

 Suppose first that the system is totally ergodic. We can assume
 that its
Kronecker factor has the form $(G,\mathcal{G},m,R_b)$ where $G$ is
a connected compact abelian group, $\mathcal{G}$ is the Borel
$\sigma$-algebra, $m$ is the Haar measure, and $R_bx=x+b$ for some
$b\in G$. If $V$ is an open subset of $G\times G$, define $
S=\{n\in\N\colon (q_1(n)b, q_2(n)b) \in V\}. $
 We claim
that if $f_1,f_2, f_3\in L^\infty(\mu)$ are such that
$\E(f_i|\mathcal{K})=0$ for some $i=1,2,3$, then
\begin{equation}\label{E:S'}
\lim_{N-M\to\infty} \frac{1}{N-M} \sum_{n=M}^{N-1}{\bf 1}_{S}(n)
\cdot T^{k_1q_1(n)}f_1\cdot T^{k_2q_2(n)}f_2\cdot
T^{l_1q_1(n)+l_2q_2(n)}f_3=0
\end{equation}
where the limit is taken in $L^2(\mu)$. We verify this as follows:
First notice that since ${\bf 1}_S(n)={\bf 1}_V(q_1(n)b,q_2(n)b)$
and the function ${\bf 1}_V$ is Riemann integrable, using an
approximation argument  it suffices to show that \eqref{E:S'}
holds for $\chi_1(q_1(n)b)\cdot \chi_2(q_2(n)b) $ in place of
${\bf 1}_S(n)$, where $\chi_1,\chi_2$ are any two characters of
$G$. To see this consider the transformation $T'=T\times
R_{b/k_1}\times R_{b/k_2}$ acting on $G^3$ with the Haar measure,
where by $b/k$ we denote a solution to the equation $kx=b$ (since
$G$ is connected such a solution always exists). Let $T'_t$,
$t\in[0,1]$, be the ergodic components of $T'$.  It is well known
(\cite{Fu1}) that  if $\E(f_1|\mathcal{K}(T))=0$ then
$\E(f_1(x)\cdot \chi_1(y)|\mathcal{K}(T'_t))=0$ for a.e. $t$.
Applying  part $(iii)$ of Theorem~B for the ergodic components of
$T'$ in place of $T$, and the functions $f_1(x)\cdot \chi_1(y)$,
$f_2\cdot \chi_2(z)$, $f_3$ in place of $f_1,f_2,f_3$, we get the
advertised identity.

We will now apply \eqref{E:S'} for the set
$$
S_\delta=\{n\in\N\colon
(q_1(n)b,q_2(n)b)\in B(0,\delta)\times B(0,\delta)\}
$$ in place of
$S$, where $\delta>0$. First notice that since the transformation
$R_b$ is totally ergodic and the polynomials $q_1,q_2$ are
linearly independent, the sequence
$\{(q_1(n)b,q_2(n)b)\}_{n\in\N}$ is uniformly distributed in
$G\times G$. Hence,
\begin{equation}\label{jnm0}
\lim_{N-M\to\infty}\frac{|S_\delta\cap
[M,N)|}{N-M}=\big(m(B(0,\delta))\big)^2.
\end{equation}
 So \eqref{E:S'} immediately gives that
 if $f_0,f_1,f_2,f_3\in
L^\infty(\mu)$ then
\begin{multline}\label{E:S2'}
\lim_{N-M\to\infty} \frac{1}{|S_\delta\cap [M,N)|} \sum_{n\in
S_\delta\cap [M,N)} \int f_0 \cdot T^{k_1q_1(n)}f_1\cdot
T^{k_2q_2(n)}f_2
\cdot T^{l_1q_1(n)+l_2q_2(n)}f_3\ d\mu=\\
\lim_{N-M\to\infty} \frac{1}{|S_\delta\cap [M,N)|} \sum_{n\in
S_\delta\cap [M,N)} \int \tilde{f}_0 \cdot
R_b^{k_1q_1(n)}\tilde{f}_1\cdot R_b^{k_2q_2(n)}\tilde{f}_2  \cdot
R_b^{l_1q_1(n)+l_2q_2(n)}\tilde{f}_3 \ dm
\end{multline}
where $\tilde{f}_i=\E(f_i|\mathcal{K})$ for $i=0,1,2,3$.

We claim that the second limit in \eqref{E:S2'} is equal to
\begin{equation}\label{E:S3'}
\frac{1}{\big(m(B(0,\delta))\big)^2}\int_{B(0,\delta)\times
B(0,\delta)} \int_G \tilde{f}_0(g)\cdot \tilde{f}_1(g+k_1t)\cdot
\tilde{f}_2(g+k_2s) \cdot \tilde{f}_3(g+l_1t+l_2s)\ dm(g)\  dm(t)
\ dm(s).
\end{equation}
One can verify this as follows:  Since the sequence
$\{(q_1(n)b,q_2(n)b)\}_{n\in\N}$ is uniformly distributed in
$G\times G$,  for every Riemann integrable function $F\colon
G\times G\to\C$ one has
\begin{equation}\label{jnm}
\lim_{N-M\to\infty}\frac{1}{N-M} \sum_{n=M}^{N-1}
F(q_1(n)b,q_2(n)b)=\int_{G\times G} F(t,s) \ dm(t)\ dm(s).
\end{equation}
 Applying
\eqref{jnm} for
$$
F(t,s)={\bf 1}_{B(0,\delta)\times B(0,\delta)}(t,s) \cdot \int_G
\tilde{f}_0(g) \cdot \tilde{f}_1(g+k_1t)\cdot
\tilde{f}_2(g+k_2s)\cdot \tilde{f}_3(g+l_1t+l_2s)\ dm(g),
$$
and using \eqref{jnm0} immediately gives that the limit in
\eqref{E:S2'} is equal to  \eqref{E:S3'}.

So we are left with estimating \eqref{E:S3'}, for some
appropriately chosen $\delta$. First notice that if $F\colon
G\times G\to \C$ is continuous then
$$
\lim_{\delta\to 0}
\frac{1}{\big(m(B(0,\delta))\big)^2}\int_{B(0,\delta)\times
B(0,\delta)} F(t,s) \ dt\ ds=F(0,0).
$$
So if $\delta$ is small enough,  and  $f_i=f={\bf 1}_A$, for
$i=0,1,2,3$,  the quantity in \eqref{E:S3'} is greater than
$$
\int (\tilde{f})^4 \ dm -\varepsilon\geq \Big(\int \tilde{f}\
dm\Big)^4-\varepsilon=\mu(A)^4-\varepsilon.
$$
Summarizing, we have shown that if $W(p_1,p_2,p_3)\leq 2$  and the
system $(X,\mathcal{B},\mu,T)$ is totally ergodic, then for every
$\varepsilon>0$, if $\delta$ is small enough we have
$$
\lim_{N-M\to\infty} \frac{1}{|S_\delta\cap [M,N)|} \sum_{n\in
S_\delta\cap [M,N)} \mu(A\cap T^{-p_1(n)}A\cap T^{-p_2(n)}A\cap
T^{-p_3(n)}A)\geq \mu(A)^4-\varepsilon.
$$
This completes the proof of Case $1$ for totally ergodic systems.

In the general case, since the Kronecker factor $\mathcal{K}$ is
an inverse limit of systems with finite rational Kronecker factor
$\mathcal{K}_{rat}$, we can choose $r\in \N$ and a factor
$\mathcal{K}'$ of $\mathcal{K}$ such that $\mathcal{K}'\cap
\mathcal{K}_{rat}=\mathcal{K}_r$ and
$$
\norm{\E(f|\mathcal{K}')-\E(f|\mathcal{K})}_{L^2(\mu)}\leq
\varepsilon/3.
$$
Then up to an error term $\varepsilon$, equation  \eqref{E:S2'}
remains valid after replacing $\tilde{f}_i$ with
$\E(f_i|\mathcal{K}')$, for $i=0,1,2,3$. The system
$(\mathcal{K}',T)$ is isomorphic to an ergodic rotation on
$H\times\Z_r$, where $H$ is a connected abelian group. We write
$q_i(rn)=rq_i'(n)$, $i=1,2$, for some integer polynomials
$q_1,q_2$, and work with $T^r$ in place of $T$ and $q_i'(n)$ in
place of $q_i(n)$, $i=1,2$. Arguing as in the totally ergodic case
we get the desired lower bound, completing the proof of Case $1$.

{\bf Case 2.} Suppose that some permutation of the essentially distinct
polynomials has the form $\{lp,mp,(l+m)p\}$ for some integer
polynomial $p$ with $p(0)=0$ and $l,m\in\Z$.
Our tactic will be  similar to the one used in the previous case
but extra complications arise because the relevant characteristic
factor in this case is not ``abelian''.

Suppose first that the system is totally ergodic.  Using an
approximation argument we can assume that the factor
$\mathcal{Z}_2$ is  an ergodic $2$-step nilsystem, say
$(X=G/\Gamma,\mathcal{G}/\Gamma,m,T_a)$. By Proposition~\ref{P:te}
we have that  $X$ is connected. Since $G$ is $2$-step nilpotent,
the subgroup $\Gamma_2=G_2\cap\Gamma$ is normal in $G$. So
$G/\Gamma_2$ is a group and $X=(G/\Gamma_2)/(\Gamma/\Gamma_2)$.
Using this representation for $X$ we can assume that
$\Gamma_2=\{e\}$ and so $G_2$ is a compact abelian Lie group.
Since $G_2$ is connected we can further assume that it is a finite
dimensional torus with the Haar measure $\lambda$. Likewise,
$Z=X/[G,G]$ is a connected compact abelian group and so we can
assume that it is a finite dimensional torus with the Haar measure
$\lambda'$.

If $\pi\colon X\to Z$ is the natural projection, and $V$ is an
open subset of $Z$, let $ S=\{n\in \N\colon p(n)a_0\in V\} $ where
$a_0=\pi(a\Gamma)$ (we use  additive notation on $Z$). We first
claim that if $f_1,f_2,f_3\in L^\infty(\mu)$ are such that
$\E(f_i|\mathcal{Z}_2)=0$ for some $i=1,2,3$, and $l,m,r$ are
distinct nonzero integers, then
\begin{equation}\label{E:S}
\lim_{N-M\to\infty} \frac{1}{N-M} \sum_{n=M}^{N-1}{\bf 1}_{S}(n)
\cdot T_a^{lp(n)}f_1\cdot T_a^{mp(n)}f_2 \cdot T_a^{rp(n)}f_3=0
\end{equation}
where the limit is taken in $L^2(\mu)$. We verify this as follows:
We can assume that the integers $l,m,r$ are relatively prime (if
not we write $l=l'd,m=m'd,r=r'd$ where $d=\text{gcd}(l,m,r)$, and
work with the polynomial family $\{l'p',m'p',rp'\}$ where $p'=dp$).
Hence,  there exist $l_1,m_1,r_1\in\Z$ such that $ll_1+mm_1+rr_1=1$.
 Since ${\bf 1}_S(n)={\bf 1}_V(p(n)a_0)$ and the function ${\bf 1}_V$ is
Riemann integrable, using an approximation argument it suffices to
verify \eqref{E:S} with $\chi(p(n)a_0)$ in place of ${\bf
1}_S(n)$, where $\chi$ is any character of $Z$ (using our notation
we have $\chi(a^{p(n)}\Gamma)=\chi(p(n)a_0)$). This last statement
follows immediately by applying Theorem~A for the functions
$f_1\cdot \chi(l_1g),f_2\cdot\chi(m_1g),f_3\cdot\chi(r_1g)$ in
place of $f_1,f_2,f_3$.

We will now apply \eqref{E:S} for the set
$$
S_\delta=\{n\in\N\colon p(n)a_0\in B(0,\delta)\}
$$
 in place of $S$, where $\delta>0$. First notice that since the sequence
$p(n)a_0$ is uniformly distributed in $Z$ we have that
\begin{equation}\label{kkl}
\lim_{N-M\to\infty}\frac{|S_\delta\cap [M,N)|}{N-M}=
\lambda'(B(0,\delta))=m(B(G_2,\delta))>0.
\end{equation}
So \eqref{E:S} immediately gives that    for $f_0,f_1,f_2,f_3\in
L^\infty(\mu)$ we have
\begin{multline}\label{E:S2''}
\lim_{N-M\to\infty} \frac{1}{|S_\delta\cap [M,N)|} \sum_{n\in
S_\delta\cap
[M,N)} \int f_0 \cdot T_a^{lp(n)}f_1\cdot T_a^{mp(n)}f_2\cdot
T_a^{rp(n)}f_3\ d\mu=\\
\lim_{N-M\to\infty} \frac{1}{|S_\delta\cap [M,N)|} \sum_{n\in
S_\delta\cap [M,N)} \int \tilde{f}_0 \cdot
T_a^{lp(n)}\tilde{f}_1\cdot T_a^{mp(n)}\tilde{f}_2\cdot
T_a^{rp(n)}\tilde{f}_3\ dm
\end{multline}
where $\tilde{f}_i=\E(f_i|\mathcal{Z}_2)$ for $i=0,1,2,3$. We
claim that the second limit in \eqref{E:S2''} is  equal
to\footnote{It may not be immediately obvious but  the next
integral is well defined. For more details see \cite{Z}.}
\begin{multline}\label{E:S3''}
\frac{1}{m(B(G_2,\delta))} \int_X \int_{B(G_2,\delta)} \int_{G_2}
\tilde{f}_0(g\Gamma)\cdot
\tilde{f}_1(gg_1^lg_2^{\binom{l}{2}}\Gamma)\cdot
\tilde{f}_2(gg_1^mg_2^{\binom{m}{2}}\Gamma) \cdot\\
\tilde{f}_3(gg_1^{r}g_2^{\binom{r}{2}}\Gamma)\ d\lambda(g_2) \
dm(g_1\Gamma) \ dm(g\Gamma),
\end{multline}
where $B(G_2,\delta)=\pi^{-1}(B(0,\delta))$.
This can be seen as follows: Since $X$ is connected, we can use
the formula of Theorem~\ref{T:Z} with $p(n)$ in place of $n$ (by
Theorem~A), $\chi, \tilde{f}_1,\tilde{f}_2,\tilde{f}_3$ in place
of $f_1,f_k,f_l,f_m$, and $1$ in place of all other $f_i$,  where
$\chi$ is any character of $Z$. We get that for a.e. $x=g\Gamma\in
X$ we have
\begin{multline}\label{qwe}
\lim_{N-M\to\infty} \frac{1}{N-M} \sum_{n=M}^{N-1}
\chi(p(n)a_0)\cdot \tilde{f}_1(a^{lp(n)}x)\cdot
\tilde{f}_2(a^{mp(n)}x)\cdot
\tilde{f}_3(a^{rp(n)}x)=\\
\int_X \int_{G_2}  \chi(g_1\Gamma)\cdot
\tilde{f}_1(gg_1^lg_2^{\binom{l}{2}}\Gamma)\cdot
\tilde{f}_2(gg_1^mg_2^{\binom{m}{2}}\Gamma) \cdot
\tilde{f}_3(gg_1^rg_2^{\binom{r}{2}}\Gamma) \ d\lambda(g_2) \
dm(g_1\Gamma).
\end{multline}
 Using an approximation argument we can verify that \eqref{qwe} holds with ${\bf
1}_{B(G_2,\delta)}$ in place of $\chi$. If we multiply this last
identity with $f_0(g\Gamma)$, then integrate with respect to
$m(x)$ and use \eqref{kkl}, we get that the limit in
\eqref{E:S2''} is equal to \eqref{E:S3''}, proving our claim.

So we are left with estimating \eqref{E:S3''} for $r=l+m$ and some
well chosen $\delta>0$. It suffices to show that
 when all functions are equal to $f={\bf 1}_A$
 the limit of \eqref{E:S3''} as $\delta\to 0$  is
greater or equal than $\mu(A)^4$.  Since
$\pi^{-1}(0)=(G_2\Gamma)/\Gamma \simeq G_2$ it is not hard to see
that this limit is equal to
$$
\int_X \int_{G_2\times G_2}  \tilde{f}(g\Gamma)\cdot
\tilde{f}(gg_1^lg_1^{\binom{l}{2}}\Gamma)\cdot
\tilde{f}(gg_1^mg_2^{\binom{m}{2}}\Gamma) \cdot
\tilde{f}(gg_1^{l+m}g_2^{\binom{l+m}{2}}\Gamma)\ d\lambda(g_2) \
d\lambda(g_1) \ dm(g\Gamma).
$$
Since the elements of $G_2$ commute with all the elements of $G$
we can write the last integral as
\begin{equation}\label{E:S5}
\int_X \int_{G_2\times G_2}  \tilde{f}(x)\cdot
\tilde{f}(g_1^lg_1^{\binom{l}{2}}x)\cdot
\tilde{f}(g_1^mg_2^{\binom{m}{2}}x) \cdot
\tilde{f}(g_1^{l+m}g_2^{\binom{l+m}{2}}x)\ d\lambda(g_1) \
d\lambda(g_2) \ dm(x).
\end{equation}
It will be more convenient  to rewrite \eqref{E:S5} as
\begin{equation}\label{E:S6}
\int_X  \int_{G_2\times G_2\times G_2} \tilde{f}(gx)\cdot
\tilde{f}(gg_1^lg_1^{\binom{l}{2}}x)\cdot
\tilde{f}(gg_1^mg_2^{\binom{m}{2}}x) \cdot
\tilde{f}(gg_1^{l+m}g_2^{\binom{l+m}{2}}x)\ d\lambda(g)\
d\lambda(g_1) \ d\lambda(g_2) \ dm(x).
\end{equation}
An easy algebraic manipulation shows that the set
$$
\{(g,gg_1^lg_2^{\binom{l}{2}},gg_1^mg_2^{\binom{m}{2}},
gg_1^{l+m}g_2^{\binom{l+m}{2}})\colon g, g_1,g_2\in G_2\}
$$
is equal to the set
$$
\{(h_1,h_2,h_3,h_4)\in G_2^4\colon
h_1^{m-l}h_3^{l+m}=h_4^{m-l}h_2^{l+m}\}.\footnote{The symmetry
of this equation is what allows us to obtain the required lower
bounds. This symmetry is violated when
 $r\neq l+m$ making
it rather unlikely that similar lower bounds hold in this case. We
discuss this more in Section~\ref{counterexamples}.}
$$
So the integral \eqref{E:S6} can be rewritten as
$$
\int_X \int_{G_2}\Big(\int_{h_1^{m-l}h_3^{m+l}=h}
\tilde{f}(h_1x)\cdot \tilde{f}(h_3x)\ d\lambda(h_1,h_3)\Big)^2
d\lambda(h)\ dm(x).
$$
Using Cauchy-Schwarz and a change of variables, we see that the
last integral is greater or equal than
$$
\int_X \Big(\int_{G_2} \tilde{f}(hx) \ d\lambda(h)\Big)^4  dm(x),
$$
which is greater or equal than
$$
\Big(\int_X \int_{G_2} \tilde{f}(hx) \ d\lambda(h)\ dm(x)\Big)^4=
\Big(\int_X \tilde{f}(x) \ dm(x)\Big)^4=\mu(A)^4.
$$
This completes the proof for totally ergodic systems.

In the general case, since every nilsystem  is an inverse limit of
nilsystems with finite rational Kronecker factor
$\mathcal{K}_{rat}$, there exists  $r_0\in \N$ and a factor
$\mathcal{Y}$ of our system, such that $\mathcal{Y}$ is a
nilsystem, $\mathcal{Y}\cap \mathcal{K}_{rat}=\mathcal{K}_{r_0}$,
and
$$
\norm{\E(f|\mathcal{Y})-\E(f|\mathcal{Z}_2)}_{L^2(\mu)}\leq
\varepsilon/10.
$$
Then up to an error term $\varepsilon$ equation  \eqref{E:S2''}
remains valid after replacing $\tilde{f}_i$ with
$\E(f_i|\mathcal{Y})$, for $i=0,1,2,3$. Moreover, the ergodic
components of the system $(\mathcal{Y},T^{r_0})$ are totally
ergodic. We write $p(r_0n)=r_0q(n)$ for some integer polynomial
$q$ and work with $T^{r_0}$ in place of $T$ and $q(n)$ in place of
$p(n)$. Arguing as in the totally ergodic case we get the desired
lower bound.
\end{proof}

\subsection{Conditional counterexamples for  the exceptional cases}\label{counterexamples}
We explain why we expect the lower bounds of Theorems~C and C' to
fail for the exceptional polynomial families $(e_1)$ with $l<m<r$
and $r\neq l+m$,  $(e_2)$, $(e_3)$. To avoid unnecessary
complications we will work out the details for two typical cases,
the general case can be treated in a similar fashion.

We first review a notion  defined by Ruzsa in \cite{R}. We consider
equations of the form
\begin{equation}\label{E:homogeneous}
 a_1x_1+\ldots +a_kx_k=0,
\end{equation}
where $a_1,\ldots,a_k\in\Z$  satisfy $a_1+\ldots+a_k=0$. Let
$\Lambda_N$ be the subset of maximum cardinality of
$\{1,\ldots,N\}$ that does not contain  solutions to
\eqref{E:homogeneous} with $k$ distinct entries. We define the
\emph{type}  of the equation \eqref{E:homogeneous} to be the
number
 $$
 \Gamma=\limsup_{N\to\infty}
 \frac{\log{|\Lambda_N|}}{\log{N}}.
 $$
For example the equation $ax+by=az+bw$ with $a,b\neq 0$ has type
$1/2$ (see \cite{R}), and an example of Behrend~\cite{Beh} shows
that the  equations $x+y=2z$, $x+y+z=3w$ have type $1$.
It seems plausible that the equation
\begin{equation}\label{type}
ax+by=cz+dw, \quad a,b,c,d>0, \quad a\neq c,d, \quad a+b=c+d,
\end{equation}
has always type $1$.
In support of this, very recently it was shown in \cite{Ko}  that
``most'' equations  of the form \eqref{type}  have type  $1$ (for
example the equation $3x+y=2z+2w$). It also seems plausible that
any equation in five variables (in fact in any odd number of
variables) with nonzero coefficients has type $1$.


\begin{proposition}\label{P:dfd}
$(i)$ Suppose that the equation $x+8z=6y+3w$ has type
$\Gamma$.\footnote{This equation is not covered by the results of
\cite{Ko}, so it is not yet known whether it has type $1$.}
Then for
every $\delta<\Gamma$ there exists an ergodic system
$(X,\mathcal{X},\mu,T)$ and $A\in\mathcal{X}$ such that
$$
\mu(A\cap T^{-2n}A \cap T^{-3n} A\cap T^{-4n} A)< \mu(A)^c
$$
for every $n\in\N$, where $c=\frac{2-\delta}{1-\delta}$ $($notice
that $c\to\infty $ as $\delta\to 1$$)$.

$(ii)$ Suppose that the equation $2x+y+w=2z+2v$ has type $\Gamma$.
Then for every $\delta<\Gamma$  there exists an ergodic system
$(X,\mathcal{X},\mu,T)$ and $A\in\mathcal{X}$ such that
$$
\mu(A\cap T^{-n}A \cap T^{-2n} A\cap T^{-n^2} A)< \mu(A)^d
$$
for every $n\in\N$, where $d=\frac{1}{2}\cdot
\frac{2-\delta}{1-\delta}$ $($notice that $d\to\infty $ as $\delta
\to 1$$)$.
\end{proposition}
\begin{proof}
$(i)$  On $X=\T^2$ with  the Borel $\sigma$-algebra $\mathcal{B}$
and the Haar measure $\mu=\lambda\times\lambda$, define the
ergodic measure preserving transformation $T\colon\T^2\to\T^2$ by
$$
T(t,s)=(t+\alpha,s+2t+\alpha),
$$
for some irrational $\alpha$. One then easily finds that
$$
T^n(t,s)=(t+n\alpha,s+2nt+n^2\alpha)
$$
for every $n\in\N$. Let $A=\T\times B\in \B$, where the set $B$
will be chosen later. We compute
\begin{multline}\label{hkl}
\mu(A\cap T^{2n}A\cap T^{3n}A\cap T^{4n}A)=\\
\int_{\T^2} {\bf 1}_B(s)\cdot {\bf 1}_B(s+4nt+4n^2\alpha) \cdot
{\bf 1}_B(s+6nt+9n^2\alpha)\cdot {\bf 1}_B(s+8nt+16n^2\alpha)  \
d\lambda(s)\  d\lambda(t) .
\end{multline}
Notice that the four numbers
$$
x=s,\ y=s+4nt+4n^2\alpha,\ z=s+6nt+9n^2\alpha,\
w=s+8nt+16n^2\alpha
$$ satisfy the equation
\begin{equation}\label{qew1}
x+8z=6y+3w.
\end{equation}
 By our assumption, for every $\delta<\Gamma$,
there exist  sets   $\Lambda_N\subset \{1,\ldots,N\}$, such that
$|\Lambda_N|\gg N^{\delta}$ and $\Lambda_N$ contains no solution
to \eqref{qew1} with distinct entries. Let
$$
B=\bigcup_{j\in\Lambda_N}
\big[\frac{j}{9N},\frac{j}{9N}+\frac{1}{81N}\big)\subset \T.
$$
Because of the condition on $\Lambda_N$ it can be easily verified
that if $x,y,z,w\in B$ satisfy \eqref{qew1} then at least two of
the $x,y,z,w$ belong to the same subinterval $I_N=
\big[\frac{j_0}{9N},\frac{j_0}{9N}+\frac{1}{81N}\big)$ for some
$j_0\in \Lambda_{N}$. Say for example that $x,y$ are these two
elements. We get  that $4nt\in
\big[-4n^2\alpha-\frac{1}{81N},-4n^2\alpha+\frac{1}{81N}\big)$ and
so $t$ belongs to a set of measure at most $2/(81N)$. The other
five cases give a similar condition, so $t$ belongs to a set
$I_{n,N}$ of measure at most $12/(81N)<1/N$.
Hence, the integral in \eqref{hkl} is at most
$$
\int_{I_{n,N}}\int_\T {\bf 1}_B(s) \ d\lambda(s)\ d\lambda(t)=
\lambda(B)\cdot \lambda(I_{n,N}) \leq \frac{|\Lambda_N|}{N^2}.
$$
Since $|\Lambda_N|\gg N^{\delta}$, for
$c=\frac{2-\delta}{1-\delta}$ an easy computation shows that
$$
\frac{|\Lambda_N|}{N^2}\ll
\Big(\frac{|\Lambda_N|}{81N}\Big)^c=\mu(A)^c.
$$
By choosing $N$ large enough we get the advertised estimate.

$(ii)$ Let $(X,\mathcal{X},\mu,T)$ be  the system used  in (i) and
let $A=B\times B\in\mathcal{B}$, where the $B$ will be chosen
later. We find
\begin{multline}\label{hkl2}
\mu(A\cap T^{n}A\cap T^{2n}A\cap T^{n^2}A)\leq\\
\int_{\T^2} {\bf 1}_B(t)\cdot {\bf 1}_B(s)\cdot {\bf
1}_B(s+2nt+n^2\alpha) \cdot {\bf 1}_B(s+4nt+4n^2\alpha)\cdot {\bf
1}_B(t+n^2\alpha) \ d\lambda(s)\  d\lambda(t) .
\end{multline}
The five numbers
$$
x=t,\ y=s, \ z=s+2nt+n^2\alpha,\ w=s+4nt+4n^2\alpha,\
v=t+n^2\alpha
$$ satisfy the equation
\begin{equation}\label{qew2}
2x+y+w=2z+2v.
\end{equation}
 By our assumption, for every $\delta<\Gamma$,
there exist sets   $\Lambda_N\subset \{1,\ldots,N\}$ such that
$|\Lambda_N|\gg N^{\delta}$ and $\Lambda_N$ contains no solution
to \eqref{qew2} with distinct entries.
By \cite{Beh}, there exist sets $B_N\subset \{1,\ldots,N\}$ with
$|B_N|\gg N e^{-C\sqrt{\log{N}}}$ and $B_N$  contains no
nontrivial $3$-term arithmetic progressions. An averaging argument
shows that there exists  $n=n(N)\in \{1,\ldots,N\}$ such that  the
set $\Lambda'_{N}=\Lambda_N\cap (B_N+n)$ has cardinality
$|\Lambda'_{N}|\geq |\Lambda_N|\cdot |B_N|$. Then $|\Lambda'_N|\gg
N^{\delta'}$ whenever  $\delta'<\delta$. Hence, for every
$\delta<\Gamma$ there exists a set $\Lambda_N\subset
\{1,\ldots,N\}$ with the following properties: $|\Lambda_N|\gg
N^{\delta}$, $\Lambda_N$ contains no solution to \eqref{qew2} with
distinct entries, and  $\Lambda_N$ contains no nontrivial $3$-term
arithmetic progressions. Let
$$
B=\bigcup_{j\in\Lambda_N}
\big[\frac{j}{4N},\frac{j}{4N}+\frac{1}{16N}\big)\subset \T.
$$
Because of the condition on $\Lambda_N$, it can be easily verified
that if $x,y,z,w,v\in B$ satisfy \eqref{qew2}, then at least two
of the $x,y,z,w,v$, excluding  the pair $x,v$, belong to the same
subinterval $[\frac{j_0}{4N},\frac{j_0}{4N}+\frac{1}{16N})$ for
some $j_0\in \Lambda_{N}$.
It follows that the integral in \eqref{hkl2} is bounded by the sum
of $\binom{5}{2}-1=9$ integrals, one of which (a typical one) has
the form
$$
I_{x,z}=\int_{A_{x,z}} {\bf 1}_B(t)\cdot {\bf 1}_B(s)\cdot {\bf
1}_B(s+2nt+n^2\alpha) \cdot {\bf 1}_B(s+4nt+4n^2\alpha)\cdot {\bf
1}_B(t+n^2\alpha) \ d\lambda(s)\  d\lambda(t),
$$
 where
$$
A_{x,z}=\big\{(t,s)\in\T^2\colon x=t, z=s+2nt+n^2\alpha\in
I_j=\big[\frac{j}{4N},\frac{j}{4N}+\frac{1}{16N}\big) \text{ for
some } j\in\Lambda_N\big\}.
$$
 We have
$$
I_{x,z}\leq \sum_{j\in\Lambda_N}  \int_{\T^2} {\bf 1}_{I_j\times
I_j}(t,s+2nt+n^2\alpha) \ d\lambda(s)\  d\lambda(t).
$$
Since the transformation $(t,s)\to (t,s+2nt+n^2\alpha)$ acting on
$\T^2$ is measure preserving, this leads to the estimate
$$
I_{x,z}\leq \sum_{j\in\Lambda_N} \lambda(I_j)^2=
\frac{|\Lambda_N|}{16^2N^2}.
$$
Similarly, we find the same bound for the other $8$ integrals.
Combining all $9$ integrals we get that
$$
\mu(A\cap T^{n}A\cap T^{2n}A\cap T^{n^2}A)\leq
\frac{|\Lambda_N|}{N^2}.
$$
Since $|\Lambda_N|\gg N^{\delta}$, for $d=\frac{1}{2}\cdot
\frac{2-\delta}{1-\delta}$  an easy computation shows that
$$
\frac{|\Lambda_N|}{N^2}\ll
\Big(\frac{|\Lambda_N|}{16N}\Big)^{2d}=\mu(A)^d.
$$
By choosing $N$ large enough we get the advertised estimate.
\end{proof}

We derive an analogous result in combinatorics:
\begin{proposition}
$(i)$ Suppose that the equation $x+8z=6y+3w$ has type $\Gamma$.
Then for every $\delta<\Gamma$  there exists $\Lambda\subset \N$
such that
$$
d^*(\Lambda\cap (\Lambda +2n) \cap (\Lambda+3n)\cap (\Lambda+4n))<
d^*(\Lambda)^c
$$
for every $n\in\N$,  where $c=\frac{2-\delta}{1-\delta}$.

$(ii)$ Suppose that the equation $2x+y+w=2z+2v$ has type $\Gamma$.
Then for every $\delta<\Gamma$ there exists $\Lambda\subset \N$
such that
$$
d^*(\Lambda\cap (\Lambda +n) \cap (\Lambda+2n)\cap (\Lambda+n^2))<
d^*(\Lambda)^d
$$
for every $n\in\N$, where $d=\frac{1}{2}\cdot
\frac{2-\delta}{1-\delta}$.
\end{proposition}
\begin{proof}
Let $(X,\mathcal{X},\mu,T)$ be the system and $A\in\mathcal{X}$ be
the set used in the proof of Proposition~\ref{P:dfd}.  Fix an
$x_0\in X$ and let
$$
\Lambda=\{n\in\N\colon T^nx_0\in A\}.
$$
Since the system is uniquely ergodic and the  boundary of the set
$A\cap T^{n_1}A\cap\cdots\cap T^{n_k}A$ has measure zero, we have
that
\begin{align*}
\mu(A\cap T^{n_1}A\cap\cdots\cap T^{n_k}A)=&
\lim_{N-M\to\infty}\frac{|\{m\in[M,N)\colon T^{m}x\in A\cap
T^{n_1}A\cap\cdots\cap T^{n_k}A\}|}{N-M}\\
=& \ d^*(\Lambda\cap (\Lambda+n_1)\cap \cdots\cap (\Lambda+n_k))
\end{align*}
for every $n_1,\ldots,n_k\in \N$. The result follows from
Proposition~\ref{P:dfd}.
\end{proof}
From the previous results we conclude that if the type of the
equation $x+8z=6y+3w$ is greater than $2/3$, or the type of the
equation $2x+y+w=2z+2v$ is greater than $6/7$, then the lower
bounds of Theorems~C and C' fail for the families $\{2n,3n,4n\}$,
$\{n,2n,n^2\}$ correspondingly. If the type of both equations is
$1$ then they fail for any fixed power of $\mu(A)$ or
$d^*(\Lambda)$.

All the other  exceptional families of Theorems~C and C' can be
treated similarly. Polynomial families of the form $(e_1)$ with
$l<m<r$ and $r\neq l+m$ lead to equations of the form
\eqref{type}, and polynomials families of the form $(e_2)$,
$(e_3)$ lead to equations in five variables. Unfortunately,
none of these equations can be treated using   the results in \cite{Ko}.


\end{document}